%% file: main.tex
\pdfoutput=1
\documentclass[toc]{mathprint}
\usepackage{rutarmacros}
\usepackage{macros}

\title[Assouad spectrum of carpets]{Assouad spectrum of Gatzouras--Lalley carpets}

\author[Banaji]
  {Amlan Banaji}
  {
      Department of Mathematics and Statistics,
      University of Jyväskylä,
      P.O.\ Box 35 (MaD),
      FI-40014 University of Jyväskylä,
      Finland
  }
  {banajimath@gmail.com}

\author[Fraser]
  {Jonathan M.\ Fraser}
  {University of St Andrews, Mathematical Institute, St Andrews, KY16 9SS, Scotland}
  {jmf32@st-andrews.ac.uk}

\author[Kolossváry]
  {István Kolossváry}
  {HUN-REN Alfréd Rényi Institute of Mathematics, 1053 Budapest, Reáltanoda u.\ 13--15, Hungary}
  {istvanko@renyi.hu}

\author[Rutar]
  {Alex Rutar}
  {
      Department of Mathematics and Statistics,
      University of Jyväskylä,
      P.O.\ Box 35 (MaD),
      FI-40014 University of Jyväskylä,
      Finland
  }
  {alex@rutar.org}

\begin{document}
\begin{abstract}
    We study the fine local scaling properties of a class of self-affine fractal sets called Gatzouras--Lalley carpets.
    More precisely, we establish a formula for the Assouad spectrum of all Gatzouras--Lalley carpets as the concave conjugate of an explicit piecewise-analytic function combined with a simple parameter change.
    Our formula implies a number of novel properties for the Assouad spectrum not previously observed for dynamically invariant sets; in particular, the Assouad spectrum can be a non-trivial differentiable function on the entire domain $(0,1)$ and can be strictly concave on open intervals.
    Our proof introduces a general framework for covering arguments using techniques developed in the context of multifractal analysis, including the method of types from large deviations theory and Lagrange duality from optimisation theory.

    \vspace{0.5cm}

    \emph{Key words and phrases.} Gatzouras--Lalley carpet, Assouad spectrum, Assouad dimension, box dimension, method of types, Lagrange duality

    \vspace{0.2cm}

    \emph{2020 Mathematics Subject Classification.} 28A80 (Primary); 37C45, 49N15 (Secondary)
\end{abstract}

\section{Introduction}\label{s:introduction}
\subsection{Overview}
The \emph{Assouad dimension} is a notion of dimension which captures the worst-case exponential scaling of a set at all locations and between all pairs of scales.
The origin of this concept can be found in the work of Assouad \cite{zbl:0396.46035} in his study of bi-Lipschitz embeddability of metric spaces in Euclidean space.
This dimension was also implicit in earlier work of Furstenberg \cite{zbl:0208.32203,zbl:1154.37322}, though in a different form as the \emph{star dimension} only recently to be shown to be equivalent \cite{doi:10.1093/imrn/rnw336}.
Especially over the past few decades, the Assouad dimension has received a large amount of attention from a variety of perspectives; we refer the reader to the books \cite{zbl:1467.28001,zbl:1222.37004,zbl:1201.30002} and the many references cited therein.

One downside of the Assouad dimension, however, is that it ignores any intermediate scaling properties of a set.
At perhaps the other extreme, the more familiar notion of the \emph{box dimension} captures the average scaling while ignoring any exceptional scaling on small parts of the set.
In order to balance these two features in a meaningful way, Fraser \& Yu introduced the \emph{Assouad spectrum} \cite{zbl:1390.28019}, which is a continuously parameterised family of dimensions bounded below by the upper box dimension and above by the Assouad dimension.
More precisely, the Assouad spectrum of a compact set $K\subset\R^d$ at $\theta\in(0,1)$ is defined by
\begin{align*}
    \dimAs\theta K = \inf\Bigl\{&\alpha:(\exists C>0)\,(\forall 0<R<1)\,(\forall x\in K)\\*
                                &N_{R^{1/\theta}}\bigl(B(x,R)\cap K\bigr)\leq C\left(\frac{R}{R^{1/\theta}}\right)^\alpha\Bigr\}.
\end{align*}
Here, for a general set $E\subset\R^d$ and $r>0$, $N_r(E)$ denotes the smallest number of balls $B(x,r)$ with $x\in E$ required to cover the set $E$.
As $\theta$ converges to $0$, $\dimAs\theta K$ converges to the upper box dimension of $K$, and as $\theta$ converges to $1$, $\dimAs\theta K$ converges to the \emph{quasi-Assouad dimension} of Lü \& Xi \cite{zbl:1345.28019}, as proven in \cite{zbl:1410.28008}.
In general, the quasi-Assouad dimension and Assouad dimension need not agree, though they do coincide for the sets under consideration in this document.

Beyond being a useful bi-Lipschitz invariant, the Assouad spectrum has also played an important role in applications: for instance, it plays a key role in the study of $L^p$-improving properties of spherical maximal operators with restricted dilation sets \cite{zbl:1561.42022,arxiv:2501.12805}, and was used to prove sharp quasiconformal distortion estimates for polynomial spirals \cite{zbl:1549.30098}.

The Assouad spectrum has been explicitly computed for a wide variety of well-studied dynamical sets with distinct box and Assouad dimensions, such as some overlapping self-similar sets and self-affine sets \cite{zbl:1407.28002,zbl:1561.28063}, Bedford--McMullen carpets \cite{zbl:1407.28002}, certain Kleinian limit sets \cite{zbl:1520.28005}, parabolic Julia sets \cite{arxiv:2203.04943}, elliptical polynomial spirals \cite{zbl:1510.28006}, and random sets such as Mandelbrot percolation \cite{zbl:1407.28002}.
In all of these examples, the Assouad spectrum has a very particular form: it is piecewise convex, with pieces of the form $\theta\mapsto a+b\cdot\frac{\theta}{1-\theta}$ for some $a,b\in\R$, and with points of non-differentiability corresponding to certain ``geometrically meaningful'' values.
The Assouad spectrum has also been computed for certain infinitely generated self-conformal sets, where more complicated forms can appear based on the fine structure of the set of fixed points (which may be prescribed in a non-dynamical way) \cite{arxiv:2207.11611}.

On the other hand, a characterisation of the possible forms of the Assouad spectrum is known: a general function $h(\theta)$ can be the Assouad spectrum of a set $E\subset\R^d$ if and only if $h$ satisfies a simple and explicit functional equation \cite{arxiv:2206.06921}.
This functional equation is very flexible and in particular a wide variety of behaviour for the Assouad spectrum is possible in general.
\begin{figure}[t]
    \centering
    \begin{subcaptionblock}{.47\textwidth}
        \centering
        \input{figures/carpet_maps/fig}
        \caption{Maps}
        \label{f:carpet-maps}
    \end{subcaptionblock}%
    \begin{subcaptionblock}{.47\textwidth}
        \centering
        \input{figures/carpet_attractor/fig}
        \caption{Attractor}\label{sf:bar}
        \label{f:carpet-attractor}
    \end{subcaptionblock}%
    \caption{Generating maps and attractor associated with a Gatzouras--Lalley IFS.}
    \label{f:carpets}
\end{figure}

In this paper we focus on the problem of computing the Assouad spectrum for a particular class of self-affine sets, namely the \emph{Gatzouras--Lalley carpets} first introduced in \cite{zbl:0757.28011}.
These generalise the well-studied Bedford--McMullen carpets, which were introduced in \cite{BedfordThesis} and \cite{zbl:0539.28003}.
The precise definition of the class of Gatzouras--Lalley carpets is deferred to \cref{ss:gl-defs}, but a depiction of the rectangles representing the generating maps and the corresponding attractor for a specific example in this class is given in \cref{f:carpets}.
Roughly speaking, the difference between Bedford--McMullen and Gatzouras--Lalley carpets is that all rectangles for a Bedford--McMullen carpet have the same widths and heights, while Gatzouras--Lalley carpets allow inhomogeneity across the widths of columns and also of the heights of rectangles within columns.
However, we emphasize that the difference between Bedford--McMullen carpets and Gatzouras--Lalley carpets is more than just technical.
For instance, it is known that Gatzouras--Lalley carpets need not have a unique measure of maximal dimension \cite{zbl:1228.37025}.
More famously, the simplest known example of a self-affine set with no ergodic measure of maximal dimension is defined using a 3-dimensional analogue of the Gatzouras--Lalley construction \cite{zbl:1387.37026}.
These features are not exhibited by Bedford--McMullen sponges in any dimension.

Our main contribution is to establish an explicit formula for the Assouad spectrum of a Gatzouras--Lalley carpet as the concave conjugate of the minimum of a finite family of elementary functions combined with a simple parameter change---see \cref{it:spectrum-formula} for the precise statement.
We observe novel and unexpected behaviour which has not been previously witnessed in any dynamically-invariant examples.
More precisely, we obtain the following direct qualitative consequences of our explicit formula:
\begin{enumerate}[nl]
    \item The Assouad spectrum can be a non-trivial differentiable function of $\theta$.
        In fact, we give a simple characterisation of differentiability, and observe that non-differentiability is generic (both topologically and measure theoretically) within the parameter space of Gatzouras--Lalley carpets---see \cref{ic:phasetrans}.
    \item The Assouad spectrum can have a non-trivial interval on which it is strictly concave.
        Again, this property holds generically---see \cref{ic:curvature}
    \item The Assouad spectrum is increasing and piecewise analytic, but with potentially arbitrarily many phase transitions and phase transitions of arbitrarily high odd order---see \cref{r:higherphase}.
\end{enumerate}
These properties are direct consequences of the explicit formula, and their proofs can be found in \cref{s:consequences}.
The proof of the explicit formula constitutes the majority of the work in this paper.

In contrast to the Bedford--McMullen case (for which the derivation of the Assouad spectrum is, generally speaking, straightforward), the Gatzouras--Lalley case has substantially more technical complications as a result of the inhomogeneity between columns and within each column.
Our proof of the general formula uses recent techniques developed in the context of multifractal analysis, which we highlight here.
First, in \cref{s:variational} we establish a general \emph{variational formula} for the Assouad spectrum as a certain constrained maximisation problem involving information-theoretic quantities evaluated at Bernoulli measures.
A key tool here is the \emph{method of types} from large deviations theory, which was used in \cite{zbl:1549.37013} to compute the $L^q$-spectrum of self-affine sponges and in \cite{arxiv:2111.05625} to calculate the intermediate dimensions of Bedford--McMullen carpets.
The explicit covering arguments build on and refine the fine covering strategies for self-affine carpets used in \cite{zbl:1278.37032,zbl:1305.28021,arxiv:2309.11971v1}.
Secondly, in \cref{s:explicit} we solve this variational formula to obtain our explicit formula.
The main complexity here is that the variational formula is a non-smooth and non-convex optimisation problem.
Our key technique here is the \emph{geometry of Lagrange duality}, which was used in \cite{arxiv:2312.08974} in order to elucidate the concave conjugate relationship apparent in the multifractal formalism.

We believe that these techniques will provide a general framework to handle coarse covering arguments in many settings beyond those considered in this paper.

\subsection{Gatzouras--Lalley carpets and dimensions}\label{ss:gl-defs}
Before we state our main results, we formally define Gatzouras--Lalley carpets and recall the known formulae for the box and Assouad dimensions.

Fix a finite index set $\mathcal{I}$ and for each $i \in \mathcal{I}$ fix constants $0<b_i<a_i<1$ and $c_i,d_i\in\R$.
Consider the diagonal self-affine iterated function system (IFS) $\{T_i\}_{i\in\mathcal{I}}$ where $T_i\colon\R^2\to\R^2$ is given by
\begin{equation}\label{e:rectangularmap}
    T_i(x,y) = (a_i x+ c_i, b_i y+d_i).
\end{equation}
By Hutchinson's application of the contraction mapping principle \cite{zbl:0598.28011}, the IFS $\{T_i\}_{i\in\mathcal{I}}$ has a unique non-empty compact \emph{attractor} $K$, which satisfies \begin{equation*}
K=\bigcup_{i\in\mathcal{I}}T_i(K).
\end{equation*}

Now, let $\eta\colon\R^2\to\R$ denote the orthogonal projection onto the first coordinate, that is $\eta(x,y)=x$.
For each $i\in\mathcal{I}$, we define the projected map $S_i$ as the unique homothety which satisfies $\eta\circ T_i= S_i\circ\eta$, or equivalently $S_i(x)=a_i x+c_i$.
Thus $\eta$ induces a map $\eta\colon\mathcal{I}\to\eta(\mathcal{I})$ taking each $i\in\mathcal{I}$ to the equivalence class $\eta(i)=\{j\in\mathcal{I}:S_j=S_i\}$, which we refer to as the \emph{column} containing index $i$.
By convention, we will refer to elements of $\mathcal{I}$ using index letters (such as $i$, $j$) and to elements of $\eta(\mathcal{I})$ using underlined index letters (such as $\ih$, $\jh$).
Note that the expressions $S_{\ih}$ and $a_{\ih}$ are well-defined for $\ih\in\eta(\mathcal{I})$.

Recall that an IFS $\{F_i\}_{i\in\mathcal{J}}$ satisfies the \emph{open set condition} with respect to an open set $U$ if $F_i(U)\subset U$ and $F_i(U)\cap F_j(U)=\varnothing$ for all $i\neq j\in\mathcal{J}$.
The following definition concerns the main class of self-affine sets studied in \cite{zbl:0757.28011}.
\begin{definition}
    We say that the IFS $\{T_i\}_{i\in\mathcal{I}}$ is \emph{Gatzouras--Lalley} if:
    \begin{enumerate}[nl,r]
        \item\label{im:osc} the original IFS $\{T_i\}_{i\in\mathcal{I}}$ satisfies the open set condition with respect to $(0,1)^2$, and
        \item\label{im:osc-proj} the projected IFS $\{S_{\ih}\}_{\ih\in\eta(\mathcal{I})}$ satisfies the open set condition with respect to $(0,1)$.
    \end{enumerate}
\end{definition}
As depicted in \cref{f:carpets}, the key features of a Gatzouras--Lalley IFS are that the rectangles $T_i([0,1]^2)$ cannot overlap except possibly along edges, they lie in columns which themselves cannot overlap except possibly along edges, and the height of each rectangle is strictly less than its width.

Next, we describe known formulae for the box and Assouad dimensions, and in the process introduce some of the notation that will be needed to state the Assouad spectrum formula.
We recall that the \emph{box dimension} of a compact set $K$ is defined to be
\begin{equation*}
    \dimB K = \lim_{r\to 0} \frac{\log N_r(K)}{\log (1/r)},
\end{equation*}
when the limit exists, and the \emph{Assouad dimension} of $F$ is defined as
\begin{align*}
    \dimA K = \inf\Bigl\{&\alpha:(\exists C>0)\,(\forall 0<r<R<1)\,(\forall x\in K)\\*
                                &N_{r}\bigl(B(x,R)\cap K\bigr)\leq C\left(\frac{R}{r}\right)^\alpha\Bigr\}.
\end{align*}
Of course, by fixing the upper scale $R$, we observe that $\dimB K\leq\dimA K$ in general.

Now, let $K$ be a \emph{Gatzouras--Lalley carpet} (i.e.\ the attractor of a Gatzouras--Lalley IFS).
First, the projection $\eta(K)$ is the attractor of the projected IFS $\{S_{\ih}\}_{\ih\in\eta(\mathcal{I})}$, which is a self-similar set satisfying the open set condition, and therefore has box dimension determined by the equation
\begin{equation*}
    \sum_{\jh\in\eta(\mathcal{I})}a_{\jh}^{\dimB\eta(K)}=1.
\end{equation*}
Next, let $t_{\min}$ denote the unique solution to
\begin{equation*}
    \sum_{\jh\in\eta(\mathcal{I})}\sum_{i\in\eta^{-1}(\jh)}a_{\jh}^{\dimB\eta(K)}b_i^{t_{\min}}=1.
\end{equation*}
We denote this quantity by $t_{\min}$ because of the usage in \cref{e:tau-formula} below.
We interpret $t_{\min}$ as the ``average'' column dimension, weighted appropriately using the column widths $a_{\jh}$.
Then the box dimension of $K$ was calculated in \cite[Theorem~2.4]{zbl:0757.28011} as
\begin{equation*}
    \dimB K=\dimB\eta(K)+t_{\min}.
\end{equation*}
Finally, for each $\jh\in\eta(\mathcal{I})$, define $s_{\jh}$ and $t_{\max}$ by the rules
\begin{equation*}
    \sum_{i\in\eta^{-1}(\jh)}b_i^{s_{\jh}}=1\qquad\text{and}\qquad t_{\max}=\max_{\jh\in\eta(\mathcal{I})}s_{\jh}.
\end{equation*}
In other words, $s_{\jh}$ is the dimension of the attractor of the IFS consisting only of the maps in column $\jh$, and $t_{\max}$ is the maximal column dimension.
Then it was proven in \cite{zbl:1278.37032} that
\begin{equation*}
    \dimA K=\dimB\eta(K)+t_{\max}.
\end{equation*}

We observe that
\begin{equation*}
    0\leq \min_{\jh\in\eta(\mathcal{I})}s_{\jh}\leq t_{\min}\leq \max_{\jh\in\eta(\mathcal{I})}s_{\jh}=t_{\max}\leq 1.
\end{equation*}
Moreover, if either the second or third inequalities are equalities, all notions of dimension for $K$ under consideration in this paper coincide (in fact, $K$ is Ahlfors--David regular).

\subsection{Main formula}
We now state our explicit formula for the Assouad spectrum of a Gatzouras--Lalley carpet.

First, for $\jh\in\eta(\mathcal{I})$ and $t\in\R$, define
\begin{equation*}
    \psi_{\jh}(t)=\frac{\log\sum_{i\in\eta^{-1}(\jh)}b_i^t}{\log a_{\jh}}.
\end{equation*}
Note that $\psi_{\jh}$ is strictly increasing and concave with unique zero $s_{\jh}$.
The Assouad spectrum of $K$ will be described in terms of the \emph{column pressure}
\begin{equation}\label{e:tau-formula}
    \tau(t)=
    \begin{cases}
        \min_{\jh\in\eta(\mathcal{I})}\psi_{\jh}(t) &: t\in [t_{\min},t_{\max}]\\*
        -\infty &:\text{ otherwise}.
    \end{cases}
\end{equation}
Observe that $\tau$ is strictly increasing on $[t_{\min},t_{\max}]$ and $\tau(t_{\max})=0$ is the unique zero of $\tau$, with $\tau(t)<0$ for $t_{\min}\leq t<t_{\max}$.
Moreover, $\tau$ is a minimum of concave functions, and is therefore concave.
We denote its concave conjugate by
\begin{equation*}
    \tau^*(\alpha)=\inf_{t\in \R}(t\alpha-\tau(t)).
\end{equation*}
Finally, define the parameter change
\begin{equation*}
    \phi(\theta)=\frac{1/\theta-1}{1-1/\kappa_{\max}}\qquad\text{where}\qquad\kappa_{\max}=\max_{i\in\mathcal{I}}\frac{\log b_i}{\log a_i}.
\end{equation*}
Our main result is the following formula for the Assouad spectrum of a Gatzouras--Lalley carpet.
\begin{itheorem}\label{it:spectrum-formula}
    Let $\{T_i\}_{i\in\mathcal{I}}$ be a Gatzouras--Lalley IFS with attractor $K$.
    Then for all $\theta\in(0,1)$,
    \begin{equation*}
        \dimAs\theta K=\dimB\eta(K)+\frac{\tau^*(\phi(\theta))}{\phi(\theta)}.
    \end{equation*}
\end{itheorem}
In \cref{s:consequences}, we derive a more transparent formula for the Assouad spectrum, and deduce many interesting qualitative features of the spectrum from it.
The proof of \cref{it:spectrum-formula} is subsequently given in \cref{s:variational} and \cref{s:explicit}.

\subsection{Notation}
We will find it useful to use various forms of asymptotic notation.
Given functions $f,g\colon A\to\R$, we write $f\lesssim g$ if there is a constant $C>0$ such that $f(a)\leq C g(a)$ for all $a\in A$.
We write $f\approx g$ if $f\lesssim g$ and $f\gtrsim g$.
We will also use Landau's $O$ notation: we say that $f=O(g)$ if there is a constant $C>0$ so that $|f(a)|\leq C|g(a)|$ for all $a\in A$.
The constants in the asymptotic notation will always be allowed to implicitly depend on the underlying IFS and a fixed parameter $\theta\in(0,1)$.
Any other dependence will be explicitly indicated by a subscript, such as $\lesssim_\varepsilon$ or $O_\varepsilon$.

\section{Consequences of the main formula}\label{s:consequences}
Throughout this section, we recall the notation introduced in \cref{s:introduction}, in particular the notation used in the statement of \cref{it:spectrum-formula}.

\subsection{An alternative formula for the Assouad spectrum}
We begin by introducing some notation to decompose the function $\tau$ in a meaningful way.

First, we distinguish a particular type of column.
\begin{definition}
    We say that a column $\jh\in\eta(\mathcal{I})$ is \emph{homogeneous} if there is a $b_{\jh}$ so that $b_i=b_{\jh}$ for all $i\in\eta^{-1}(\jh)$.
\end{definition}
For example, when $K$ is a Bedford--McMullen carpet, every column is homogeneous.

Homogeneity is characterised by the following elementary lemma, the proof of which follows directly from the definition of $\psi_{\jh}$.
\begin{lemma}\label{l:homog}
    A column $\jh\in\eta(\mathcal{I})$ is homogeneous if and only if $\psi_{\jh}$ is affine.
    Moreover:
    \begin{enumerate}[nl,r]
        \item If $\psi_{\jh}$ is affine, then
            \begin{equation*}
                \psi_{\jh}(t)=t\cdot\frac{\log b_{\jh}}{\log a_{\jh}}+\frac{\log \#\eta^{-1}(\jh)}{\log a_{\jh}} = \kappa_{\jh} \cdot (t-s_{\jh}).
            \end{equation*}
        \item If $\psi_{\jh}$ is not affine, then it is strictly concave.
    \end{enumerate}
\end{lemma}

Now write
\begin{equation}\label{e:gt-def-intro}
    g(t)=\min_{\jh\in\eta(\mathcal{I})}\psi_{\jh}(t)
\end{equation}
which is a minimum of analytic functions (note the close relationship with the function $\tau$).
By analyticity, for all $\ih,\jh\in\eta(\mathcal{I})$, either $\psi_{\ih}=\psi_{\jh}$ or the set $\{t\in\R:\psi_{\ih}(t)=\psi_{\jh}(t)\}$ is finite.
Thus there exists a partition
\begin{equation*}
    t_{\min}=t_0<t_1<\cdots<t_m=t_{\max}
\end{equation*}
of the interval $[t_{\min},t_{\max}]$, with corresponding parts $I_n=[t_{n-1},t_{n}]$, such that for each $n\in\N$, there is a $\jh_n$ so that
\begin{equation*}
    g(t)=\begin{cases}
        \psi_{\jh_1}(t) &: t\in I_1\\*
                        \quad\vdots\\*
        \psi_{\jh_m}(t) &: t\in I_m
    \end{cases}
\end{equation*}
and moreover for all $1\leq n\leq m-1$, $\psi_{\jh_n}\neq\psi_{\jh_{n+1}}$.
The latter property ensures that the partition $(I_n)_{n=1}^m$ is uniquely determined.
We refer to this partition as the \emph{spectrum partition} associated with the IFS $\{T_i\}_{i\in\mathcal{I}}$.
We associate with the spectrum partition the following additional information, all of which depends only on the underlying IFS.
\begin{enumerate}[nl]
    \item We say that a part $I_n$ is \emph{homogeneous} if $\jh_n$ is homogeneous, and \emph{inhomogeneous} otherwise.
    \item We associate with each part $I_n$ the function $g_n=\psi_{\jh_n}$, which is analytic on the open interval $I_n^\circ$.
    \item We associate with each part $I_n=[t_{n-1},t_{n}]$ the endpoint derivatives
        \begin{equation}\label{e:thetaminmaxdef}
            \theta_{n,\min}=\phi^{-1}(g_n'(t_{n-1}))\qquad\text{and}\qquad\theta_{n,\max}=\phi^{-1}(g_n'(t_n)).
        \end{equation}
\end{enumerate}
Expanding the definitions of $\phi^{-1}$ and $g$, we may equivalently write
\begin{equation*}
    \theta_{n,\min}=\frac{1}{\partial^+g(t_{n-1})\cdot(1-1/\kappa_{\max})+1}\quad\text{and}\quad
    \theta_{n,\max}=\frac{1}{\partial^-g(t_{n})\cdot(1-1/\kappa_{\max})+1}.
\end{equation*}
Here, and elsewhere, for a concave function $f$, we write $\partial^-f(x)$ and $\partial^+f(x)$ to denote the left and right derivatives of $f$ at $x$ respectively (when they exist), and set $\partial f(x)=[\partial^+f(x),\partial^- f(x)]$.
In particular,
\begin{equation*}
    \theta_{1,\min}\leq\theta_{1,\max}\leq\theta_{2,\min}\leq\cdots\leq\theta_{n,\max}
\end{equation*}
and moreover $\theta_{n,\min}=\theta_{n,\max}$ if and only if $I_n$ is homogeneous.
We define $\theta_{\min}=\theta_{1,\min}$ and $\theta_{\max}=\theta_{m,\max}$.

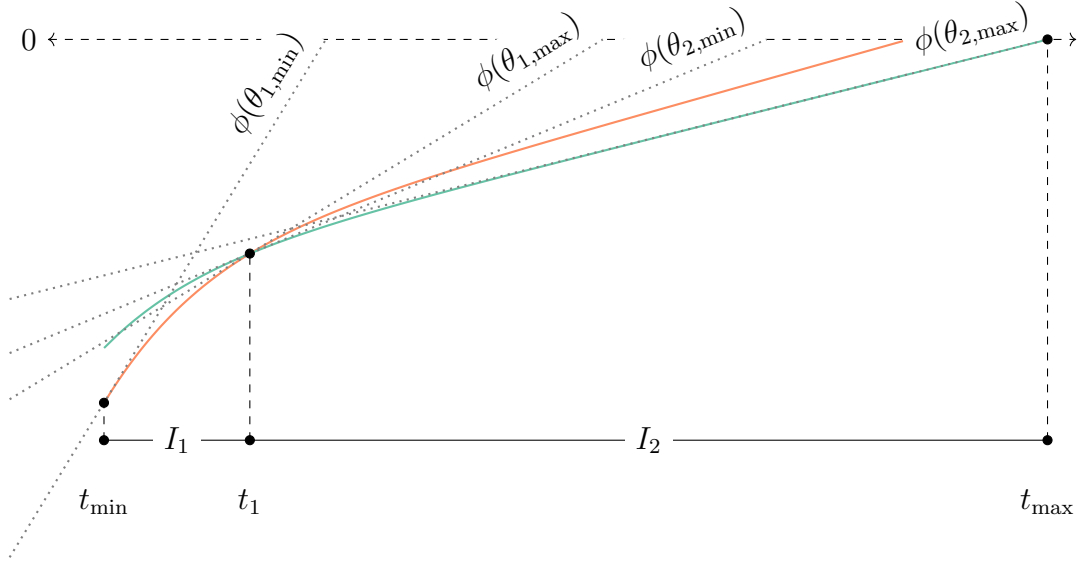
\begin{figure}[t]
    \centering
    \input{figures/piecewise/fig}
    \caption{A depiction of the spectrum partition.
        The dotted lines are tangents to the function $g=\min\{\psi_{\jh_1},\psi_{\jh_2}\}$ at the points $t_{\min}$, $t_1$, and $t_{\max}$ corresponding to the left and right derivatives, where appropriate.
        The labels indicate the slopes of the dotted lines.
    }
    \label{f:piecewise}
\end{figure}

Using this notation, \cref{it:spectrum-formula} yields the following explicit piecewise formula for the Assouad spectrum of $K$.
For clarity in the below formula, note that $g$ is strictly increasing with $g(t_{\max})=0$.
The cases of this formula are depicted in \cref{f:piecewise}.
\begin{icorollary}\label{ic:piecewise}
    Fix a Gatzouras--Lalley IFS $\{T_i\}_{i\in\mathcal{I}}$ with spectrum partition $(I_n)_{n=1}^m$, and associated data $g_n$, $\theta_{n,\min}$, and $\theta_{n,\max}$ as above.
    If $\dimB K=\dimA K$, then $K$ is Ahlfors--David regular and $\dimB K=\dimA K=\dimAs\theta K$ for all $\theta\in(0,1)$.

    Otherwise, one of the following conditions holds for each $\theta\in(0,1)$:
    \begin{enumerate}[nl,r]
        \item\label{im:box} We have $\theta\leq\theta_{\min}$.
            Then
            \begin{equation*}
                \dimAs\theta K =\dimB K-\frac{g(t_{\min})}{\phi(\theta)}.
            \end{equation*}
        \item\label{im:diff-vals} There is an $n\in\{1,\ldots,m\}$ so that $\theta_{n,\min}<\theta<\theta_{n,\max}$.
            Then
            \begin{equation*}
                \dimAs\theta K = \dimB\eta(K)+\frac{g_n^*(\phi(\theta))}{\phi(\theta)}.
            \end{equation*}
        \item\label{im:non-diff-vals} There is an $n\in\{1,\ldots,m-1\}$ so that $\theta_{n,\max}\leq\theta\leq\theta_{n+1,\min}$.
            Then
            \begin{equation*}
                \dimAs\theta K = \dimB\eta(K)+t_n-\frac{g_n(t_n)}{\phi(\theta)}.
            \end{equation*}
        \item\label{im:Assouad} We have $\theta\geq\theta_{\max}$.
            Then
            \begin{equation*}
                \dimAs\theta K=\dimA K.
            \end{equation*}
    \end{enumerate}
    Moreover,
    \begin{equation*}
        \theta_{\max}=\inf\{\theta\in(0,1):\dimAs\theta K=\dimA K\},
    \end{equation*}
    and $\theta_{\min}=\theta_{\max}$ if and only if $m=1$ and $I_1$ is homogeneous.
\end{icorollary}
\begin{proof}
    We recall the definition of $\tau$ and the main result proven in \cref{it:spectrum-formula}.
    Note that the formulae for $\theta\leq\theta_{\min}$ and $\theta\geq\theta_{\max}$ follow directly by the definition of the concave conjugate applied at the endpoint.
    Similarly, for $\theta_{\min}<\theta<\theta_{\max}$, the formulae follow directly since
    \begin{equation*}
        \dimAs\theta K=\dimB\eta(K)+\frac{g^*(\phi(\theta))}{\phi(\theta)}.
    \end{equation*}
    The parts in case \cref{im:diff-vals} correspond to the interiors of $I_n$, in which case $g^*$ and $g_n^*$ agree, and the parts in case \cref{im:non-diff-vals} correspond to the points of non-differentiability of $g$, which can only occur on the endpoints between adjacent $I_n$.
    In the latter case, $g^*$ is an affine function of $\alpha$ with an explicit formula depending only on the value of $g_n$ at $t_n$.

    Moreover, suppose $\dimB K\neq\dimA K$, so that $t_{\min}<t_{\max}$.
    First, since we recall that $s_{\jh}$ is the unique zero of $\psi_{\jh}$ for all $\jh\in\eta(\mathcal{I})$ and $t_{\max}=\max s_{\jh}$, it follows that $t_{\max}$ is the unique zero of $g$.
    Thus, for $\theta<\theta_{\max}$, let $t$ be such that $\phi(\theta)\in\partial g(t)$, so $t<t_{\max}$ and since $g$ is concave,
    \begin{equation*}
        g^*(\phi(\theta))=g(t_{\max})+g^*(\phi(\theta)) < t_{\max}\phi(\theta).
    \end{equation*}
    Dividing through by $\phi(\theta)$ gives the claim.

    It is also clear directly from the definition that $\theta_{\min}=\theta_{\max}$ if and only if $g$ is a affine function on the interval $[t_{\min},t_{\max}]$, which occurs if and only if $m=1$ and $I_1$ is homogeneous.
\end{proof}
We note that the formulae in \cref{im:box} and \cref{im:Assouad} are in fact special cases of \cref{im:non-diff-vals} after substituting the respective value of $t_n$.
In these three cases, expanding the definition of $\phi$, the formula is of the form $a+b\frac{\theta}{1-\theta}$, which has occurred previously as discussed in the introduction.
In contrast, case \cref{im:diff-vals} is novel and very much \emph{not} of this form.
This case occurs only in the presence of an inhomogeneous column.

In case \cref{im:diff-vals}, an implicit formula for the concave conjugate can be obtained from \cref{p:minimisers} (also see the discussion in \cref{r:optim}).
The spectrum can also be obtained parametrically as a function of $t$: given a part $I_n=[t_{n-1},t_n]$ with corresponding column function $\psi=\psi_{\jh_n}=g_n$, the graph of the function $\theta\mapsto \dimAs\theta K$ on $(\theta_{n,\min},\theta_{n,\max})$ is the same as the image of the interval $(t_{n-1},t_n)$ under the map
\begin{equation}\label{e:param}
    t\mapsto \left(\phi^{-1}(\psi'(t)),\dimB\eta(K)+t-\frac{\psi(t)}{\psi'(t)}\right).
\end{equation}

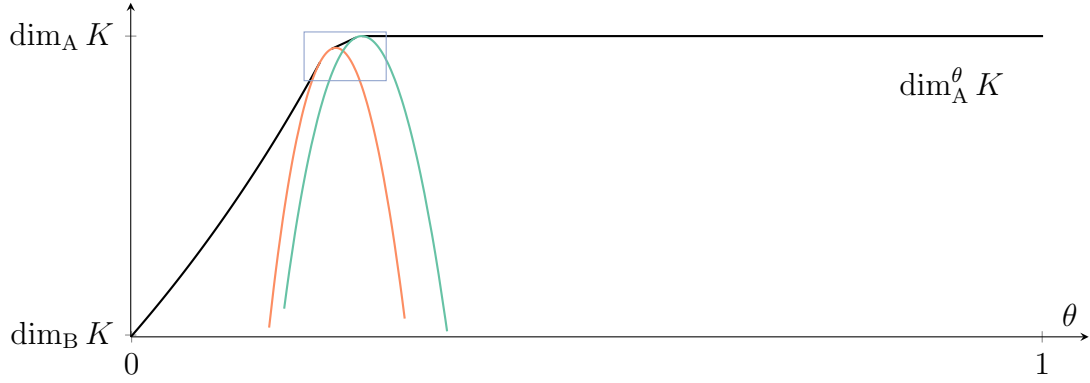
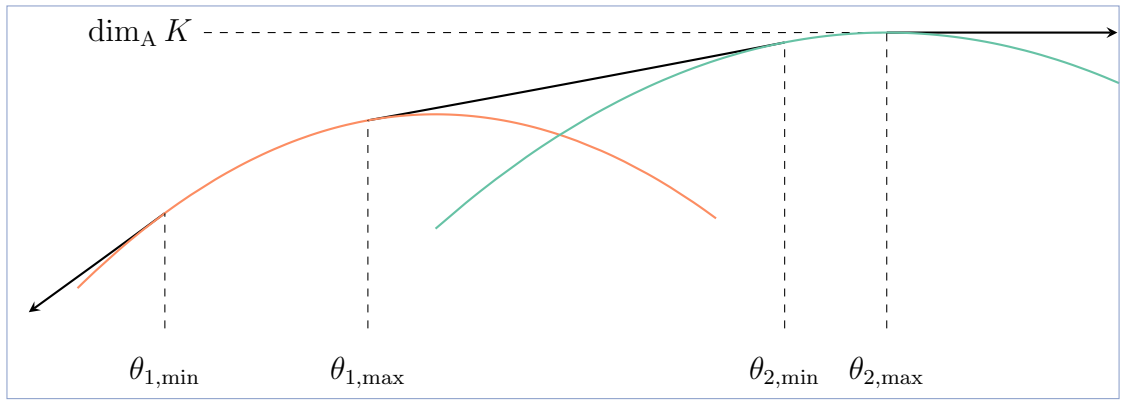
\begin{figure}[t]
    \centering
    \begin{subfigure}{\textwidth}
        \centering
        \input{figures/two_bump_spectrum/fig}
        \caption{Plot of the original spectrum.}
        \label{f:two-bump-transform}
    \end{subfigure}%
    \vspace{0.5cm}
    \begin{subfigure}{\textwidth}
        \centering
        \input{figures/two_bump_spectrum_zoomed/fig}
        \caption{Plot restricted to the rectangular region.}
        \label{f:two-bump-spectrum}
    \end{subfigure}%
    \caption{A depiction of decomposition provided by \cref{ic:piecewise}.
        The coloured curves are of the form $\dimB\eta(K)+g_i^*(\phi(\theta))/\phi(\theta)$ for $i=1,2$.
        The spectrum is differentiable but not twice differentiable at each $\theta_{i,\min}$ and $\theta_{i,\max}$ for $i=1,2$.
    }
    \label{f:spectrum-decomp}
\end{figure}

A graphical depiction of the curve $g(t)$ is given in \cref{f:piecewise}, and the decomposition provided by \cref{ic:piecewise} is given in \cref{f:spectrum-decomp}.
The case $\theta_{\min}=\theta_{\max}$ is satisfied, for example, by every Bedford--McMullen carpet (or more generally by any Gatzouras--Lalley carpet with $\frac{\log b_i}{\log a_i}$ constant for $i\in\mathcal{I}$).

Using this decomposition, we can establish that $\dimAs\theta K$ is an increasing function of $K$.
Here, we recall that the \emph{upper Assouad spectrum} is defined by
\begin{align*}
    \dimuAs\theta K = \inf\Bigl\{&\alpha:(\exists C>0)\,(\forall 0<r\leq R^{1/\theta}<R<1)\,(\forall x\in K)\\*
                                &N_{r}\bigl(B(x,R)\cap K\bigr)\leq C\left(\frac{R}{r}\right)^\alpha\Bigr\}.
\end{align*}
The main result of \cite{zbl:1410.28008} implies that if $\dimAs\theta K$ is an increasing function, then $\dimAs\theta K=\dimuAs\theta K$ for all $\theta\in(0,1)$.
\begin{icorollary}\label{ic:incr}
    Let $K$ be a Gatzouras--Lalley carpet.
    Then $\dimAs\theta K$ is strictly increasing on $(0,\theta_{\max}]$ and constant on $[\theta_{\max},1)$.
    In particular, $\dimuAs\theta K=\dimAs\theta K$ for all $\theta\in(0,1)$.
\end{icorollary}
\begin{proof}
    Using the formula in \cref{ic:piecewise}, it is clear that $\dimAs\theta K$ is strictly increasing on intervals corresponding to cases \cref{im:box} and \cref{im:non-diff-vals}.
    Moreover, $\dimAs\theta K$ is constant on \cref{im:Assouad} (which corresponds to the case when $\theta\in[\theta_{\max},1)$).
    For the remaining case \cref{im:diff-vals}, fix some $n=1,\ldots,m$, write $\psi=g_n$ and let $s$ denote the unique zero of $\psi$.
    Now let $\theta\in(\theta_{n,\min},\theta_{n,\max})$ be arbitrary and let $\alpha=\phi(\theta)$.
    Note that $\alpha>0$ and moreover $\alpha>\psi'(s)$ since $\psi$ is strictly concave (for otherwise $\theta_{n,\min}=\theta_{n,\max}$) with $\psi(s)=0$.

    Note that $\psi$ is strictly increasing and strictly concave, so $-\psi^{-1}$ is a well-defined strictly concave function.
    For notational simplicity, let
    \begin{equation*}
        \varphi=(-\psi^{-1})^*\qquad\text{and}\qquad F(\alpha)=\frac{\psi^*(\alpha)}{\alpha}.
    \end{equation*}
    In particular, by definition of the concave conjugate, using the change of variable $y=\psi(t)$,
    \begin{equation*}
        F(\alpha) = \inf_{t\in\R}\left\{t-\frac{\psi(t)}{\alpha}\right\}= \inf_{y\in\R}\left\{y\cdot\left(-\frac{1}{\alpha}\right)-(-\psi^{-1})(y)\right\}=\varphi\left(-\frac{1}{\alpha}\right).
    \end{equation*}
    Moreover, $\varphi$ is maximised at
    \begin{equation*}
        (-\psi^{-1})'(0)=-\frac{1}{\psi'(\psi^{-1}(0))}=-\frac{1}{\psi'(s)}<-\frac{1}{\alpha}.
    \end{equation*}
    Thus since $\varphi$ is strictly concave, its derivative is strictly decreasing, so
    \begin{equation*}
        \alpha^2 F'(\alpha)=\varphi'\left(-\frac{1}{\alpha}\right)<\varphi'\left(-\frac{1}{\psi'(s)}\right)=0.
    \end{equation*}
    Thus by the chain rule, since $\phi'(\theta)<0$ for all $\theta\in(0,1)$,
    \begin{equation*}
        \frac{\dx{}}{\dx{\theta}}\dimAs\theta K=F'\bigl(\phi(\theta)\bigr)\cdot\phi'(\theta)>0.
    \end{equation*}
    In particular, the spectrum is a strictly increasing function of $\theta$ on $(0,\theta_{\max}]$.
\end{proof}

\subsection{Differentiability and higher-order phase transitions}
We can also obtain precise information concerning differentiability, as well as higher-order phase transitions.
We say that a continuous real-valued function $h$ has a \emph{phase transition of $k$\textsuperscript{th} order at $\theta$} when $k\geq 1$ is the smallest integer such that the $k$\textsuperscript{th}-derivative $h^{(k)}(\theta)$  does not exist.
Recall that we denote the number of parts in the spectrum partition by $m$.
\begin{icorollary}\label{ic:phasetrans}
    The function $\theta\mapsto\dimAs\theta K$ is piecewise analytic and the set of points where the function is not analytic is given precisely by
    \begin{equation*}
        H=\{\theta:\theta=\theta_{n,\min}\text{ or }\theta=\theta_{n,\max}\quad\text{for some}\quad n=1,\ldots,m\}.
    \end{equation*}
    At each $\theta\in H$ there is a phase transition that either has odd order or order $2$.
    Moreover:
    \begin{enumerate}[nl,r]
        \item\label{im:order-1} The set of $\theta\in(0,1)$ at which $\dimAs\theta K$ has a 1\textsuperscript{st} order phase transition is given by
            \begin{equation*}
                H_1\coloneqq\{\theta:\theta=\theta_{n,\min}=\theta_{n,\max}\quad\text{for some}\quad n=1,\ldots,m\}.
            \end{equation*}
            This implies that $\dimAs\theta K$ has precisely $k$ points of non-differentiability, where $k$ is the number of $i=1,\ldots,m$ such that $I_i$ is homogeneous, and $k\leq\#\eta(\mathcal{I})-1$.
            In particular, $\dimAs\theta K$ is differentiable if and only if each $I_n$ is inhomogeneous.
        \item\label{im:odd-order} The set of $\theta\in(0,1)$ at which $\dimAs\theta K$ has a $k$\textsuperscript{th} order phase transition for some odd integer $k\geq 3$ is given by
            \begin{equation*}
                H_{\mathrm{higher}}\coloneqq \{\theta:\theta=\theta_{n,\max}=\theta_{n+1,\min}\quad\text{for some}\quad n=1,\ldots,m-1\}\setminus H_1.
            \end{equation*}
        \item\label{im:order-2} The set of $\theta\in(0,1)$ at which $\dimAs\theta K$ has a 2\textsuperscript{nd} order phase transition is given by
            \begin{equation*}
                H_2\coloneqq H\setminus (H_1\cup H_{\mathrm{higher}}).
            \end{equation*}
    \end{enumerate}
    In particular, at $\theta_{\min}$ and $\theta_{\max}$, $\dimAs\theta K$ either has a 1\textsuperscript{st} or 2\textsuperscript{nd} order phase transition.
\end{icorollary}
\begin{proof}
    Piecewise analyticity follows directly from the piecewise formula given in \cref{ic:piecewise}.
    Moreover, since the distinct parts correspond to distinct analytic curves and any intersection of distinct analytic curves must have a phase transition of some order, we obtain the formula for $H$.

    First, to see \cref{im:order-1}, by standard properties of the concave conjugate, the derivative of $\tau^*$ at $\alpha$ exists if and only if $\tau$ is strictly concave at all $t$ for which $\alpha\in\partial\tau(t)$.
    But $\tau$ fails to be strictly concave at $t$ if and only if $t\in I_n^\circ$ for a homogeneous part $I_n$, in which case $\theta=\theta_{n,\min}=\theta_{n,\max}$, as claimed.
    Since the curves $\psi_{\jh}$ are affine for a homogeneous column $\jh$ and otherwise strictly concave, and moreover there must be at least one column $\ih$ with $s_{\ih}\leq t_{\min}$, there are at most $\#\eta(\mathcal{I})-1$ points of non-differentiability.

    Next, suppose $\theta=\theta_{n,\max}=\theta_{n+1,\min}$ for some $n=1,\ldots,m-1$.
    Equivalently, $\theta=g_n'(t_n)=g_{n+1}'(t_n)$, and since $g_n-g_{n+1}$ changes sign at $t_n$, $g_n-g_{n+1}$ has a saddle point at $t_n$.
    Therefore if $k\in\N$ is minimal so that $g_n^{(k)}(t_n)\neq g_{n+1}^{(k)}(t_n)$, then $k$ an odd integer which is at least $3$.
    In particular, if $\theta\notin H_1$, then $\dimAs\theta K$ is differentiable at $\theta$, and therefore has a phase transition of odd order $k\geq 3$.

    Otherwise, suppose $\theta\in H_2$ and $\theta=\theta_{n,\max}<\theta_{n+1,\min}$.
    Since $(g^*)''=0$ on $(\phi^{-1}(\theta_{n+1,\min}),\phi^{-1}(\theta_{n,\max}))$ and $(g^*)''$ is uniformly bounded away from $0$ on $(\phi^{-1}(\theta_{n,\max}),\phi^{-1}(\theta_{n,\min}))$, and since the parameter change $g^*(\phi(\theta))/\phi(\theta)$ is smooth, the second derivative does not exist at $\theta_{n,\max}$.
    The case $\theta_{n,\max}<\theta_{n+1,\min}=\theta$ is analogous.
    Otherwise, $\theta=\theta_{\max}$, but again $\dimAs\theta K$ is constant on $[\theta_{\max},1)$ so the same argument yields non-existence of the second derivative at $\theta_{\max}$.

    Combining these two observations yields \cref{im:odd-order} and \cref{im:order-2}, and in particular that every phase transition has either odd order or order $2$.
\end{proof}

\begin{remark}\label{r:higherphase}
    It is straightforward to see that phase transitions of order $1$ and order $2$ already occur in many of the examples given later in this document.
    We sketch a construction giving a phase transition of arbitrary odd order $k\geq 3$.

    We consider a Gatzouras--Lalley carpet with three columns each of width $1/3$, the first of which contains only one map (with arbitrary height less than $1/3$).
    Let $N\in\N$ be large and let $B\subset\R^N$ denote the set of all $N$-tuples $(b_1,\ldots,b_N)$ with $0<b_i<1/3$ and $b_1+\cdots+b_N< 1$.
    Note that to any pair $(\beta,\tilde{\beta})\in B\times B$ there is an attractor $K=K(\beta,\beta')$ where the second column has contraction ratios given by $\beta$, and the third given by $\tilde\beta$.
    Fix $\beta=\tilde{\beta}=\beta_0$ for some arbitrary initial choice of $\beta_0$ whose entries are not all equal, let $t_{\min}$, $t_{\max}$ be the values corresponding to the carpet $K(\beta_0,\beta_0)$, and let $t$ satisfy
    \begin{equation*}
        t_{\min}<t<t_{\max}.
    \end{equation*}
    Exponentiating, we see that for any given $d \geq 1$, the functions $\psi$ and $\tilde{\psi}$ corresponding to the non-trivial columns have the same $d$\textsuperscript{th} derivative at $t$ if and only if
    \begin{equation*}
        \sum_{i=1}^N \bigl((\log b_i)^d b_i^t - (\log \tilde{b}_i)^d (\tilde{b}_i)^t\bigr) = 0.
    \end{equation*}
    Moreover, by the implicit function theorem, for typical choices of $(\beta_0,\beta_0)$ and for sufficiently large $N$, the set of parameters in $B\times B$ for which at least the first $k$ derivatives match at $t$ is a non-trivial submanifold of $B\times B$ containing the point $(\beta_0,\beta_0)$.
    But given that the first $k$ derivatives match, the parameters for which the first $k+1$ derivatives match is a proper submanifold, so there must exist parameters $(\beta,\tilde\beta)$ arbitrarily close to $(\beta_0,\beta_0)$ such that precisely the first $k$ derivatives of $\psi$ and $\tilde{\psi}$ agree, but the next derivative does not.

    Since the parameters $t_{\min}$ and $t_{\max}$ are continuous functions of $(\beta,\tilde{\beta})$, and since the sign of $\psi-\tilde{\psi}$ must change at $t$ since $k$ is odd, it follows that the function $g$ corresponding to the carpet $K(\beta,\tilde{\beta})$ must have a phase transition of order $k$ at $t$.
    Finally, since $k\geq 3$, $g^*$ is differentiable so the derivative of $g^*$ is the inverse of $g'$.
    Since the reparameterisation in terms of $\phi$ is smooth, it follows that $\dimAs\theta K(\beta,\tilde{\beta})$ has a phase transition of order $k$ at $\phi^{-1}(g'(t))$.

    If one instead chooses points $t_{\min}<t_1<\cdots<t_n<t_{\max}$, a similar argument gives arbitrarily many phase transitions of arbitrary odd orders at least $3$.
\end{remark}

\subsection{Convexity and concavity}
As the final result of this section, we obtain some information concerning convexity and concavity.
\begin{icorollary}\label{ic:curvature}
    The spectrum $\dimAs\theta K$ is:
    \begin{enumerate}[nl,r]
        \item\label{im:convex} Strictly convex on each interval $(\theta_{n,\max},\theta_{n+1,\min})$ for $n=1,\ldots,m-1$ as well as the interval $(0,\theta_{\min})$;
        \item\label{im:concave} Strictly concave on the interval $(\theta_{\max}-\delta,\theta_{\max})$ for some $\delta>0$ if and only if $I_m$ is inhomogeneous; and
        \item\label{im:constant} Constant on the interval $[\theta_{\max},1)$.
    \end{enumerate}
    In particular, if $\dimAs\theta K$ is not constant, then $\dimAs\theta K$ contains a non-trivial interval of convexity, and if every column is inhomogeneous, then $\dimAs\theta K$ also contains a non-trivial interval of concavity.
\end{icorollary}
\begin{proof}
    Cases \cref{im:convex} and \cref{im:constant} follow directly from the piecewise formula for the Assouad spectrum given in \cref{ic:piecewise}.
    The remaining case which requires checking is \cref{im:concave}.
    If $I_m$ is homogeneous, $\dimAs\theta K$ is strictly convex in a neighbourhood to the left of $\theta_{\max}$ by \cref{im:convex}.

    Otherwise, assume that $I_m$ is inhomogeneous with corresponding column $\jh\in\eta(\mathcal{I})$.
    Writing $\psi=\psi_{\jh}$ and continuing the computation from the proof of \cref{ic:incr} with the same definitions of $F$ and $\varphi$, for all $\theta\in I_m^\circ$,
    \begin{equation*}
        \frac{\dx{}^2}{\dx{\theta^2}}\dimAs\theta K=F''(\phi(\theta))\cdot(\phi'(\theta))^2+F'(\phi(\theta))\cdot\phi''(\theta).
    \end{equation*}
    But $F'(\phi(\theta_{\max}))=0$ and
    \begin{align*}
        F''(\phi(\theta_{\max})) &=\varphi''\left(-\frac{1}{\phi(\theta_{\max})}\right)\frac{1}{\phi(\theta_{\max})^4}-2\varphi'\left(-\frac{1}{\phi(\theta_{\max})}\right)\frac{1}{\phi(\theta_{\max})^3}<0,
    \end{align*}
    since we recall that $\varphi$ is strictly concave so $\varphi''<0$, and $\varphi'(-1/\psi'(t_{\max}))=0$ where $\psi'(t_{\max})=\phi(\theta_{\max})$.
    Thus, by continuity of the second derivative,
    \begin{equation*}
        \frac{\dx{}^2}{\dx{\theta^2}}\dimAs\theta K<0
    \end{equation*}
    for some $\delta>0$ and $\theta\in(\theta_{\max}-\delta,\theta_{\max})$, as claimed.
\end{proof}
\begin{remark}
    For most of the examples we present in this document, the Assouad spectrum is strictly concave on each non-trivial interval $(\theta_{n,\min},\theta_{n,\max})$.
    However, in \cref{ss:convex-bump}, we construct an example with a non-trivial inhomogeneous column and a non-trivial open sub-interval of $(\theta_{n,\min},\theta_{n,\max})$ on which $\dimAs\theta K$ is strictly convex.
\end{remark}

\subsection{Examples}\label{s:examples}

\subsubsection{Homogeneous carpets}
\begin{figure}[t]
    \centering
    \begin{subfigure}{\textwidth}
        \centering
        \input{figures/homog_transform/fig}
        \caption{Plot of the function $g(t)$.}
    \end{subfigure}%
    \vspace{0.5cm}
    \begin{subfigure}{\textwidth}
        \centering
        \input{figures/homog_spectrum/fig}
        \caption{Plot of the Assouad spectrum.}
    \end{subfigure}
    \caption{Plot of the Assouad spectrum corresponding to a system with 4 homogeneous columns, 3 of which are non-trivial.
        The function $g(t)$ is a minimum of affine lines, so the corresponding spectrum $\dimAs\theta K$ is a piecewise convex function.
        The slope of $g_i(t)$ corresponds to the value $\theta_i$, for $i=1,2,3$.
        The dotted lines correspond to the concave conjugates at each $t_i$ for $i=0,\ldots,3$, extended beyond the range given by the corresponding affine lines.
    }
    \label{f:homog}
\end{figure}
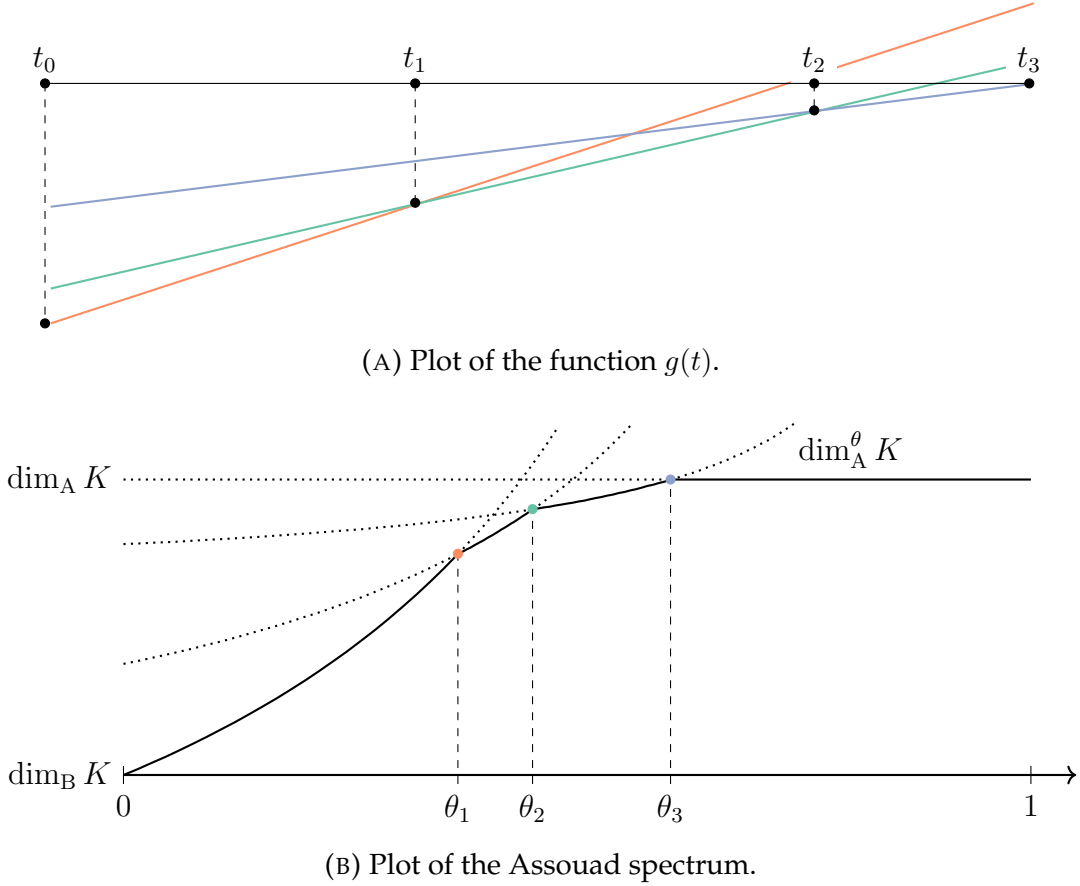
Consider the special case of Gatzouras--Lalley carpets whose columns are all homogeneous, i.e.\ for each $\jh\in \eta(\mathcal{I})$, there is a unique $b_{\jh}$ so that $b_i = b_{\jh}$ for all $i\in \eta^{-1}(\jh)$.
Recall from \cref{l:homog} that $\psi_{\jh}$ is the affine function
\begin{equation*}
    \psi_{\jh}(t) = \kappa_{\jh} \cdot (t - s_{\jh})\qquad\text{where}\qquad\kappa_{\jh}\coloneqq\frac{\log b_{\jh}}{\log a_{\jh}}
\end{equation*}
and we recall that $s_{\jh}$ is the dimension of column $\jh$.
Hence, the function $g(t)$ is piecewise affine, so for any $\theta \in(0,1)$, $\phi(\theta) \in \partial g(t_n)$ for some $n=0,\ldots,m$.
Equivalently, $\theta_n\coloneqq\theta_{n,\min}=\theta_{n,\max}$ for all $n=1,\ldots,m$, so case \cref{im:diff-vals} of \cref{ic:piecewise} never occurs.

We can also derive an explicit formula for case \cref{im:non-diff-vals} of \cref{ic:piecewise}.
Since $\psi_{\jh}$ are affine for all $\jh\in\eta(\mathcal{I})$, two functions $\psi_{\ih}$ and $\psi_{\jh}$ either have the same slope (in which case one of them does not appear in the formula at all) or intersect at
\begin{equation*}
    t_{\ih,\jh} \coloneqq \frac{\kappa_{\ih} s_{\ih}-\kappa_{\jh} s_{\jh}}{\kappa_{\ih} -\kappa_{\jh}}\qquad\text{with value}\qquad\psi_{\ih}(t_{\ih,\jh}) = \frac{s_{\ih}-s_{\jh}}{1/\kappa_{\jh} - 1/\kappa_{\ih}}.
\end{equation*}
Given the spectrum partition of the carpet, for every $n=1,\ldots, m-1$, we have $t_n=t_{\jh_n,\jh_{n+1}}$ and we can express $\theta_n$ as
\begin{equation*}
\theta_n = \phi^{-1}(\kappa_{\jh_n}) = \frac{1}{1+\kappa_{\jh_n}\cdot (1-1/\kappa_{\max})}.
\end{equation*}
To conclude, we obtain the formula
\begin{equation*}
    \dimAs\theta K= \dimB \eta(K)+ \frac{\kappa_{\jh_n} s_{\jh_n}-\kappa_{\jh_{n+1}} s_{\jh_{n+1}}}{\kappa_{\jh_n} -\kappa_{\jh_{n+1}}} + \frac{\theta}{1-\theta} (1-1/\kappa_{\max}) \frac{s_{\jh_{n+1}} - s_{\jh_n}}{1/\kappa_{\jh_{n+1}} - 1/\kappa_{\jh_n}}.
\end{equation*}
In particular, $\dimAs{\theta_n} K= \dimB \eta(K)+ s_{\jh_n}$.
The spectrum $\dimAs\theta K$ is piecewise convex with a
point of non-differentiability at each $\theta_n$.
The number of such points is bounded from above by the number of columns minus one, and this bound is clearly optimal.

A plot of the function $g$ and the Assouad spectrum for a system with four columns (three of which are non-trivial) and three phase transitions is given in \cref{f:homog}.

\subsubsection{An example with three columns and six phase transitions}\label{ss:two-int}
\begin{figure}[t]
    \centering
    \input{figures/two_int/fig}
    \caption{Spectrum plot restricted to a small domain with a given column being relevant on multiple intervals.
        Note that the values $\theta_{1,\min}$ and $\theta_{3,\max}$ are not in the domain of this image.
    }
    \label{f:two-int}
\end{figure}
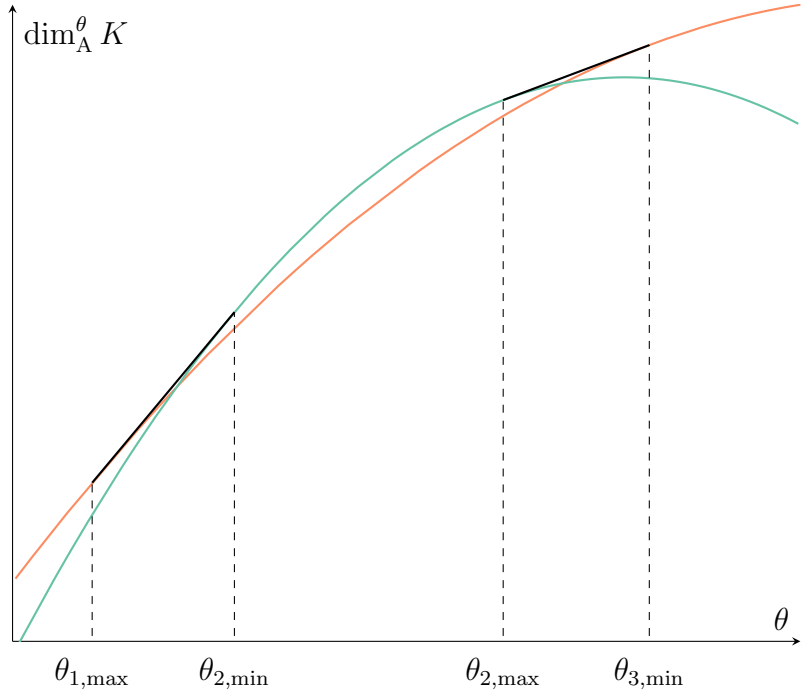
In this section, we provide an explicit example of a Gatzouras--Lalley system with 6 phase transitions, each of order 2.
Define maps
\begin{align*}
    T_{1,1}(x,y) &= (0.1\cdot x, 0.05\cdot y)\\
    T_{2,1}(x,y) &= (0.4\cdot x + 0.2, 0.00001\cdot y) & T_{2,2}(x,y) &= (0.4\cdot x + 0.2, 0.39\cdot y+0.61)\\
    T_{3,1}(x,y) &= (0.31\cdot x + 0.69, 0.000177\cdot y)& T_{3,1}(x,y) &= (0.31\cdot x + 0.69, 0.2\cdot y+0.8)
\end{align*}
This is a system with three columns, consisting of a single map $T_{1,1}$ in the first column, maps $T_{2,1}$ and $T_{2,2}$ in the second, and $T_{3,1}$ and $T_{3,2}$ in the third.
This system has the following properties, which can be determined by a straightforward (albeit tedious!) computation:
\begin{enumerate}[nl]
    \item The spectrum partition has three parts $I_1$, $I_2$, and $I_3$ where the second column dominates on the parts $I_1$ and $I_3$ and the third column dominates on part $I_2$.
    \item The Assouad spectrum has six phase transitions.
\end{enumerate}
A plot of the Assouad spectrum on a restricted domain is given in \cref{f:two-int}.

\subsubsection{Inhomogeneous column with corresponding part convex}\label{ss:convex-bump}
The Assouad spectrum corresponding to case \cref{im:diff-vals} of \cref{ic:piecewise} is often concave; in particular, this is the case for the other examples given above.
However, this is not necessarily the case in general.
To give an explicit example, we consider an IFS consisting of two columns.
The first column consists of the single map $T_{1,1}(x,y)=(4x/5, y/1000)$.
The second column has maps
\begin{equation*}
    T_{2,j}(x,y)=(x/5, b_j y)+(4/5, t_j)\qquad\text{for}\qquad j=1,\ldots,52
\end{equation*}
where $b_1=b_2=19/100$ and $b_3=\cdots=b_{52}=10^{-20}$, and the $t_j$ are chosen so that the attractor is a Gatzouras--Lalley carpet.

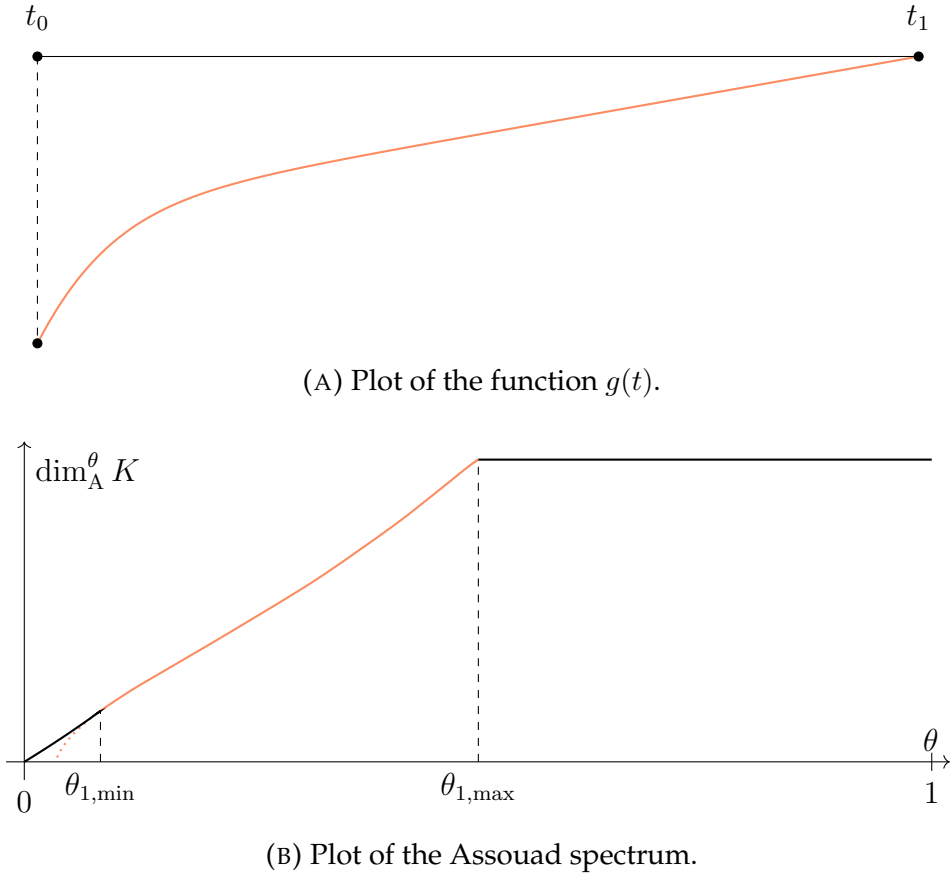
\begin{figure}[t]
    \centering
    \begin{subfigure}{\textwidth}
        \centering
        \input{figures/convex_transform/fig}
        \caption{Plot of the function $g(t)$.}
    \end{subfigure}%
    \vspace{0.5cm}
    \begin{subfigure}{\textwidth}
        \centering
        \input{figures/convex_spectrum/fig}
        \caption{Plot of the Assouad spectrum.}
    \end{subfigure}%
    \caption{Plot of the Assouad spectrum which has a convex part in the interval $(\theta_{1,\min},\theta_{1,\max})$.
        Note that the spectrum is differentiable on $(0,1)$, including at $\theta_{1,\max}$.
    }
    \label{f:convex}
\end{figure}

The convexity of the Assouad spectrum on a non-trivial open sub-interval of $(\theta_{1,\min},\theta_{1,\max})$ can be easily (albeit tediously) verified using the parametric formula \cref{e:param}.
A plot of the Assouad spectrum is given in \cref{f:convex}.

\subsection{Further work}
We now briefly describe some possible directions for future research.

\subsubsection{Natural interpretation of the concave conjugate}
The concave conjugate of the Assouad spectrum has also appeared naturally in another setting, namely local smoothing estimates with restricted times \cite{arxiv:2501.12805}.
There, the dual of the Assouad spectrum is referred to as the \emph{Legendre--Assouad function}.
It would be worthwhile to determine the relationship between this function and the column pressure $\tau(t)$ from \cref{e:tau-formula}.

\subsubsection{Bounds on number of phase transitions}
We see from \cref{ic:phasetrans} and \cref{s:examples} that it is much easier to bound the number of phase transitions when all columns are homogeneous than in the general case, and we ask the following question.
\begin{question}
For each integer $N \geq 3$, let $f(N)$ denote the supremum of $n \in \N$ for which there exists a Gatzouras--Lalley carpet with $N$ maps whose Assouad spectrum is non-analytic at $n$ distinct points.
Is $f(N)$ finite for all $N$, and if so can it be calculated or estimated?
\end{question}

\subsubsection{Assouad spectrum of more general self-affine sets}
It is natural to ask about the Assouad spectrum of more general self-affine sets.
If the condition $b_i < a_i$ is weakened to $b_i \leq a_i$ in \cref{e:rectangularmap} for some (but not all) $i$ then we expect that the main results of this paper will still hold.
More generally, it would be worthwhile to calculate the formula for carpets of the form considered by Barański \cite{zbl:1116.28008}, or for self-affine sponges in higher dimensions.

\subsubsection{The lower spectrum}
One might also compute the \emph{lower spectrum} (see \cite[§1]{zbl:1390.28019}) of a Gatzouras--Lalley carpet, which is the natural dual to the Assouad spectrum describing scaling properties of the set at locations where it is sparse.
We expect that similar techniques as used in this paper would be useful to compute a formula for the lower spectrum.

\subsubsection{Bi-Lipschitz equivalence}
Since the Assouad spectrum is a bi-Lipschitz invariant, our results immediately imply some consequences for bi-Lipschitz equivalence of self-affine carpets.
For instance, if the spectrum partition of one carpet contains a homogeneous column while the other does not, then the two carpets cannot be bi-Lipschitz equivalent.
In particular, this distinguishes Bedford--McMullen carpets from Gatzouras--Lalley carpets in which every column is inhomogeneous.

However, the Assouad spectrum does not take into account any columns which do not appear in the spectrum partition, whereas such columns may impact bi-Lipschitz equivalence.
Here, the lower spectrum might be useful, since it would also take into account the small columns.
More generally, classifying Gatzouras--Lalley carpets up to bi-Lipschitz equivalence is a challenging problem (some partial results are given in \cite{doi:10.1088/1402-4896/ad86f2}), and it would be interesting to investigate what more can be said in this direction using the Assouad and lower spectra.

\subsubsection{Intermediate dimensions of Gatzouras--Lalley carpets}
Finally, the \emph{intermediate dimensions} are another family of dimensions which depend on a parameter $\theta \in (0,1)$, and they lie between Hausdorff and box dimension.
A precise formula for the intermediate dimensions of Bedford--McMullen carpets is given in \cite{arxiv:2111.05625}, with a highly technical proof.
It is likely an even more technically challenging problem to calculate a formula for the intermediate dimensions of all Gatzouras--Lalley carpets.
However, in light of \cref{ic:phasetrans} \cref{im:order-1}, it would be interesting to know whether the intermediate dimension function is differentiable for Gatzouras--Lalley carpets whose columns are all inhomogeneous.
The intermediate dimensions of Bedford--McMullen carpets with distinct Hausdorff and box dimensions have countably many points of non-differentiability.

\section{A variational formula for the Assouad spectrum}\label{s:variational}
For the remainder of this document, we fix a Gatzouras--Lalley IFS $\{T_i\}_{i\in\mathcal{I}}$ with attractor $K$.
In this section, we establish a variational formula for the Assouad spectrum.
After setting up notation in \cref{ss:notation}, in \cref{ss:var} we state the variational formula and give a sketch of the proof.
The remaining sections are devoted to the proof of the variational formula.
We will use the \emph{method of types} from large deviations theory; an introduction can be found in \cite[§2.1.1]{zbl:1177.60035}.

\subsection{Symbolic notation and approximate squares}\label{ss:notation}
\subsubsection{Symbolic notation and coding}
Set $\mathcal{I}^*=\bigcup_{n=0}^\infty\mathcal{I}^n$ equipped with the operation of concatenation.
Given $\mtt{i}\in\mathcal{I}^n$, we denote the \emph{length} of $\mtt{i}$ by $|\mtt{i}|=n$.
We say that a word $\mtt{i}\in\mathcal{I}^*$ is a \emph{prefix} of $\mtt{k}$ if there is a word $\mtt{j}$ so that $\mtt{k}=\mtt{i}\mtt{j}$, and we write $\mtt{i}\preccurlyeq\mtt{k}$.
If $|\mtt{i}|\geq 1$, we let $\mtt{i}^-$ denote the unique prefix of $\mtt{i}$ of length $|\mtt{i}|-1$.
We denote the unique word of length 0 by $\varnothing$.

Given $\mtt{i}=(i_1,\ldots,i_n)\in\mathcal{I}^n$, we write
\begin{align*}
    a_{\mtt{i}}&=a_{i_1}\cdots a_{i_n},& b_{\mtt{i}}&=b_{i_1}\cdots b_{i_n},\\
    T_{\mtt{i}}&=T_{i_1}\circ\cdots\circ T_{i_n},&S_{\mtt{i}}&=S_{i_1}\circ\cdots\circ S_{i_n}.
\end{align*}

Next, we let $\Omega=\mathcal{I}^{\N}$ denote the space of infinite sequences on $\mathcal{I}$.
The concatenation $\mtt{i}\gamma$ for $\mtt{i}\in\mathcal{I}^*$ and $\gamma\in\Omega$ is defined similarly, so we also speak of finite prefixes of infinite words.
Given $\mtt{i}\in\mathcal{I}^*$, we denote the \emph{cylinder set}
\begin{equation*}
    [\mtt{i}]=\{\gamma\in\Omega:\mtt{i}\preccurlyeq \gamma\}.
\end{equation*}

Of course, the map $\eta\colon\mathcal{I}\to\eta(\mathcal{I})$ also induces a map on $\mathcal{I}^*$ and $\Omega$.
Sometimes, we will abuse notation and simply write $\mtt{i}\umtt{j}$ for $\mtt{i}\in\mathcal{I}^*$ and $\umtt{j}\in\eta(\mathcal{I}^*)$, in place of $\eta(\mtt{i})\umtt{j}$.
For instance, the expressions
\begin{equation*}
    a_{\mtt{i}\umtt{j}}\qquad\text{and} \qquad S_{\mtt{i}\umtt{j}}
\end{equation*}
are well-defined since they depend only on the value of $\eta(\mtt{i})$ and $\umtt{j}$.

Finally, we define the surjective coding map $\pi\colon\Omega\to K$ for $\gamma=(i_n)_{n=1}^\infty$ by
\begin{equation*}
    \{\pi(\gamma)\}=\bigcap_{n=1}^\infty T_{i_1}\circ\cdots\circ T_{i_n}(K).
\end{equation*}
As a result of the separation assumptions on a Gatzouras--Lalley IFS, we will often abuse notation and refer to sets $E\subset\Omega$ and their image $\pi(E)\subset K$ interchangeably.

\subsubsection{Notation for probability vectors}
In order to state our variational formula, we must first introduce some notation to handle the space of Bernoulli measures associated with the IFS $\{T_i\}_{i\in\mathcal{I}}$.
Let
\begin{equation*}
    \mathcal{P}=\mathcal{P}(\mathcal{I})=\Big\{(p_i)_{i\in\mathcal{I}}:p_i\geq 0,\sum_{i\in\mathcal{I}} p_i=1\Big\}\subset\R^{\mathcal{I}},
\end{equation*}
which is a compact metric space with the metric inherited from ambient Euclidean space.
We refer to probability vectors in $\mathcal{P}$ using bold-face letters, such as $\bm{w}$.
Recall that $\eta$ denotes the orthogonal projection onto the first coordinate axis.
In a similar way as before, $\eta$ induces a map $\eta\colon\mathcal{P}\to\eta(\mathcal{P})$ by the rule
\begin{equation*}
    \eta(\bm{w})=\Bigl(\sum_{i\in\eta^{-1}(\jh)}w_i\Bigr)_{\jh\in\eta(\mathcal{I})}.
\end{equation*}
Given a probability vector $\bm{w}=(w_i)_{i\in\mathcal{I}}\in\mathcal{P}$, we define the \emph{entropy}
\begin{equation*}
    H(\bm{w})=\sum_{i\in\mathcal{I}}w_i\log(1/w_i)
\end{equation*}
and \emph{Lyapunov exponents}
\begin{equation*}
    \chi_1(\bm{w})=\sum_{i\in\mathcal{I}}w_i\log(1/a_i)\qquad\text{and}\qquad \chi_2(\bm{w})=\sum_{i\in\mathcal{I}}w_i\log(1/b_i).
\end{equation*}
For $\bm{p}\in\eta(\mathcal{I})$, we define entropy in the same way, and also note that $\chi_1(\bm{p})$ is well-defined since $\chi_1(\bm{w})=\chi_1(\eta(\bm{w}))$.
Note that $H$, $\chi_1$, and $\chi_2$ are continuous positive functions on $\mathcal{P}$, and moreover $\chi_1$ and $\chi_2$ are uniformly bounded away from 0.
Since rectangles defining the Gatzouras--Lalley carpet are wider than they are tall, $\chi_1(\bm{w})<\chi_2(\bm{w})$.

Finally, we denote the \emph{logarithmic eccentricity} of $\bm{w}\in\mathcal{P}$ by
\begin{equation*}
    \Gamma(\bm{w})=\frac{\chi_2(\bm{w})}{\chi_1(\bm{w})}.
\end{equation*}
Since $0<b_i<a_i<1$ for all $i\in\mathcal{I}$, $\Gamma$ takes values in a compact interval
\begin{equation*}
    \Gamma(\mathcal{P})=[\kappa_{\min},\kappa_{\max}]\subset(1,\infty),
\end{equation*}
where
\begin{equation*}
    \kappa_{\min}=\min_{i\in\mathcal{I}}\frac{\log b_i}{\log a_i}\qquad\text{and}\qquad\kappa_{\max}=\max_{i\in\mathcal{I}}\frac{\log b_i}{\log a_i}.
\end{equation*}
We will use $\Gamma$ to measure the exponential distortion of the rectangle $T_{\mtt{i}}([0,1]^2)$ in terms of the digit frequencies corresponding to the word $\mtt{i}\in\mathcal{I}^*$, with $\kappa_{\min}$ corresponding to the cases when the rectangle looks as close as possible to a square.

\subsubsection{Pseudo-cylinders and approximate squares}
In order to perform our covering arguments, we require some symbolic notation to handle the fact that the cylinders are exponentially distorted.
A key difficulty in understanding the geometry of self-affine sets is that the cylinder sets are not squares, but instead exponentially distorted rectangles.
On the other hand, because of the vertical alignment of the columns, it turns out that we can group vertically-aligned cylinder sets together in ways which let us cover multiple cylinders simultaneously.

Our main symbolic construction is the notion of a \defn{pseudo-cylinder}.
Suppose $\mtt{i}\in\mathcal{I}^n$ and $\umtt{j}\in\eta(\mathcal{I})^m$.
Then the corresponding \defn{pseudo-cylinder} is the set
\begin{equation*}
    P(\mtt{i},\umtt{j})=\{\gamma=(i_k)_{k=1}^\infty\in\Omega:(i_1,\ldots,i_n)=\mtt{i}\text{ and }\eta(i_{n+1},\ldots,i_{n+m})=\umtt{j}\}.
\end{equation*}
Note that the map $(\mtt{i},\umtt{j})\mapsto P(\mtt{i},\umtt{j})$ is injective.
Equivalently, the pseudo-cylinder $P(\mtt{i},\umtt{j})$ is the union of the cylinders contained in $[\mtt{i}]$ which all lie in column $\eta(\mtt{i})\umtt{j}$; that is,
\begin{equation}\label{e:pseudo-cylinder-rep}
    P(\mtt{i},\umtt{j})=\bigcup_{\mtt{k}\in\eta^{-1}(\umtt{j})}[\mtt{i}\mtt{k}].
\end{equation}

Now given an infinite word $\gamma\in\Omega$, let $L_n(\gamma)$ be the minimal integer so that
\begin{equation*}
    a_{\gamma_1}\cdots a_{\gamma_{L_n(\gamma)}}< b_{\gamma_1}\cdots b_{\gamma_n}.
\end{equation*}
In other words, $L_n(\gamma)$ is chosen so that the level $L_n(\gamma)$ rectangle has approximately the same width as the height of the level $n$ rectangle.
Note that $L_n(\gamma) \geq n$, and write the first $L_n(\gamma)$ entries of $\gamma$ as $\gamma\npre{L_n(\gamma)}=\mtt{i}\mtt{j}$ where $\mtt{i}\in\mathcal{I}^n$.
We then define the \defn{approximate square} $Q_n(\gamma)\subset\Omega$ by
\begin{equation*}
    Q_n(\gamma)=P(\mtt{i},\eta(\mtt{j})).
\end{equation*}
We use the term `approximate square' because the ratio of its height to its width is bounded away from $0$ and $\infty$, independently of $\gamma$.
While different $\gamma$ may define the same approximate square, the choice of $\mtt{i}$ and $\eta(\mtt{j})$ are unique.

Now for fixed $\mtt{i}$, let $\mathcal{U}(\mtt{i})\subset\eta(\mathcal{I}^*)$ denote the set of $\umtt{j}$ so that $P(\mtt{i},\umtt{j})$ is an approximate square.
We call $\mathcal{U}(\mtt{i})$ a \emph{complete section}: for every infinite word $\underline{\gamma}\in\eta(\Omega)$, there is exactly one $\umtt{j}\in\mathcal{U}(\mtt{i})$ which is a prefix of $\underline{\gamma}$.
In particular, the approximate squares $P(\mtt{i},\umtt{j})$ are disjoint in symbolic space for fixed $\mtt{i}$.

We say that a pseudo-cylinder $P(\mtt{i},\umtt{j})$ is \emph{wide} if $\umtt{j}\preccurlyeq\umtt{k}$ for some $\umtt{k}\in\mathcal{U}(\mtt{i})$; in other words, $P(\mtt{i},\umtt{j})$ contains approximate squares of the form $P(\mtt{i},\umtt{k})$.
One can think of the wide pseudo-cylinders as ``interpolating'' between the cylinder $P(\mtt{i},\varnothing)=[\mtt{i}]$ and the approximate square $P(\mtt{i},\eta(\mtt{j}))=Q_n(\gamma)$.
If $P(\mtt{i},\umtt{j})$ is wide, then $|\umtt{j}|\lesssim |\mtt{i}|$, and moreover if $P(\mtt{i},\umtt{j})$ is an approximate square, then $|\umtt{j}|\approx|\mtt{i}|$.

Finally, we group the approximate squares based on diameter (or more precisely, height) by defining, for $0<r<1$,
\begin{equation*}
    \mathcal{S}(r)=\left\{Q=P(\mtt{i},\umtt{j}):Q\text{ is an approximate square and }b_{\mtt{i}}<r\leq b_{\mtt{i}^-}\right\}.
\end{equation*}
Note that the elements of $\mathcal{S}(r)$ are pairwise disjoint in symbolic space since the condition $b_{\mtt{i}}<r\leq b_{\mtt{i}}^-$ defines a complete section on $\mathcal{I}^*$, and $\mathcal{U}(\mtt{i})$ is a complete section for each $\mtt{i}$.

\subsection{Statement of the variational formula and proof strategy}\label{ss:var}
We can now explicitly state our variational formula.
First, define a parameter change for $\bm{v}\in\mathcal{P}$ by
\begin{equation*}
    \phi(\theta,\bm{v})=\frac{1/\theta-1}{1-1/\Gamma(\bm{v})}\qquad\text{and}\qquad\phi(\theta)=\inf_{\bm{v}\in\mathcal{P}}\phi(\theta,\bm{v})=\frac{1/\theta-1}{1-1/\kappa_{\max}}.
\end{equation*}
Note that $\phi$ is strictly decreasing in $\theta$ and $\Gamma(\bm{v})$.
We let
\begin{equation}\label{e:definedeltathin}
    \begin{split}
        \Delta_{\mathrm{thin}}(\theta)&=\{(\bm{v},\bm{w})\in\mathcal{P}\times\mathcal{P}:\phi(\theta,\bm{v})\leq\Gamma(\bm{w})\},\\*
        \Delta_{\mathrm{thick}}(\theta)&=\{(\bm{v},\bm{w})\in\mathcal{P}\times\mathcal{P}:\phi(\theta,\bm{v})\geq\Gamma(\bm{w})\}.
    \end{split}
\end{equation}
Recall that $t_{\min}=\dimB K-\dimB\eta(K)$, and set
\begin{align*}
    f_{\mathrm{thin}}(\theta,\bm{v},\bm{w})&=\dimB\eta(K)+\frac{H(\bm{w})-H(\eta(\bm{w}))}{\chi_2(\bm{w})},\\
    f_{\mathrm{thick}}(\theta,\bm{v},\bm{w})&=\dimB K+\frac{1}{\phi(\theta,\bm{v})}\left(\frac{H(\bm{w})-H(\eta(\bm{w}))-t_{\min}\chi_2(\bm{w})}{\chi_1(\bm{w})}\right).
\end{align*}
Finally, write
\begin{equation}\label{f:objective}
    f(\theta,\bm{v},\bm{w})=\begin{cases}
        f_{\mathrm{thin}}(\theta,\bm{v},\bm{w})&:(\bm{v},\bm{w})\in\Delta_{\mathrm{thin}}(\theta),\\*
        f_{\mathrm{thick}}(\theta,\bm{v},\bm{w})&:(\bm{v},\bm{w})\in\Delta_{\mathrm{thick}}(\theta).
    \end{cases}
\end{equation}
It is straightforward to check that $f_{\mathrm{thin}}=f_{\mathrm{thick}}$ on $\Delta_{\mathrm{thin}}(\theta)\cap\Delta_{\mathrm{thick}}(\theta)$, so $f$ is indeed well-defined and continuous.
\begin{theorem}\label{t:gl-variational}
    Let $K$ be a Gatzouras--Lalley carpet.
    Then for all $\theta\in(0,1)$,
    \begin{equation*}
        \dimAs\theta K=\max_{(\bm{v},\bm{w})\in\mathcal{P}\times\mathcal{P}}f(\theta,\bm{v},\bm{w}).
    \end{equation*}
\end{theorem}

We now summarise the main idea of the proof.
In order to compute $\dimAs\theta K$, it suffices to consider approximate squares.
Consider an approximate square of width $R$, say, and note that it is composed of cylinders, all of which have various heights.
For our application there are two cases:
either the cylinder is \emph{thin}, i.e.\ it has height at most $R^{1/\theta}$, or the cylinder is \emph{thick}, i.e.\ it has height at least $R^{1/\theta}$.
This corresponds to the two cases of the definition in \cref{f:objective}, and the covering strategy for each case is different.
\begin{enumerate}[nl]
    \item In the \emph{thin} case (which is handled in \cref{ss:thin}), we can simply group cylinders together (forming a \emph{pseudo-cylinder}) until the heights are approximately $R^{1/\theta}$.
        We cover such cylinders simultaneously, and the count depends on $\dimB\eta(K)$ (the $\dimB K$ count does not appear).
    \item In the \emph{thick} case (which is handled in \cref{ss:thick}), we must cover each cylinder at a scale $R^{1/\theta}$, which is smaller than the height of the cylinder.
        We note that pseudo-cylinders of height $R^{1/\theta}$ and a certain width can be realised as images of approximate squares in $K$, so we can count the number of such pseudo-cylinders in terms of $\dimB K$, and then cover each one at scale $R^{1/\theta}$ using $\dimB\eta(K)$.
\end{enumerate}
The key observation is that the covering strategy for the cylinder depends only on the digit frequency (or \emph{type}) of certain indices corresponding to each case.
In the thick case, these digit frequencies are precisely a pair $(\bm{v},\bm{w})\in\mathcal{P}\times\mathcal{P}$, where $\bm{v}$ is the digit frequency of the cylinder defining the approximate square, and $\bm{w}$ is the digit frequency of the smaller cylinder.
In the thin case, $\bm{v}$ is the same but now $\bm{w}$ corresponds instead to the pseudo-cylinder.
The precise definitions of these types are given in \cref{ss:types}.
The resolutions at which these bounds become relevant depend on the logarithmic eccentricity of both the original cylinder $\bm{v}$ and the composing cylinders $\bm{w}$.
But now the crucial observation (based on the same strategy underlying the main results in \cite{zbl:1549.37013}) which allows us to complete this argument is a combination of the following two facts:
\begin{enumerate}[nl,resume]
    \item Since the scales $R^{1/\theta}$ and $R$ are exponentially separated, the set of possible types is much smaller than the number of cylinders with each type; and
    \item The covering strategy within each type class is the same, with cost precisely corresponding to the functions $f_{\mathrm{thin}}$ and $f_{\mathrm{thick}}$.
\end{enumerate}
This means that the cost to cover the approximate square is dominated by the maximal type, yielding the variational formula for the Assouad spectrum.
In fact, we will later see that when $\theta_{\min} \leq \theta \leq \theta_{\max}$, there will be a maximising vector on $\Delta_{\mathrm{thin}}(\theta)\cap\Delta_{\mathrm{thick}}(\theta)$; the corresponding cylinders will have height approximately $R^{1/\theta}$.

\subsection{Sections for approximate squares and types}\label{ss:types}
In this section, we make the notion of a \emph{type} rigorous, which is the starting point for our covering strategy.

\subsubsection{Defining the section \texorpdfstring{$\mathcal{B}^\theta(\mtt{i},\umtt{j})$}{Bθ(i,j)}}
First, we say that a subset $\mathcal{A}\subset\mathcal{I}^*$ is a \emph{section} if for any $\gamma\in\Omega$, there is \emph{at most one} $\mtt{i}\in\mathcal{A}$ which is a prefix of $\gamma$.
Note that complete sections, as introduced earlier, are sections to which no new elements can be added.
The partial order of prefixes extends to a partial order on the family of sections, where we write $\mathcal{A}\preccurlyeq\mathcal{B}$ if for all $\mtt{i}\in\mathcal{A}$, there exists $\mtt{j}\in\mathcal{B}$ such that $\mtt{i}$ is a prefix of $\mtt{j}$.
This partial order has a \emph{meet}: that is, given a finite family of sections $\mathcal{A}_1,\ldots,\mathcal{A}_n$, there is a unique section $\mathcal{A}_1\wedge\cdots\wedge\mathcal{A}_n$ which is maximal with respect to the partial order such that
\begin{equation*}
    \mathcal{A}_1\wedge\cdots\wedge\mathcal{A}_n\preccurlyeq\mathcal{A}_i
\end{equation*}
for all $i=1,\ldots,n$.

Now fix $\theta\in(0,1)$ and an approximate square $Q=P(\mtt{i},\umtt{j})$ where $\mtt{i}\in\mathcal{I}^*$ and $\mtt{j}\in\eta(\mathcal{I}^*)$.
Recall that $\eta^{-1}(\umtt{j})\subset\mathcal{I}^*$ is in bijection with the set of cylinders composing the approximate square $Q$ (see \cref{e:pseudo-cylinder-rep}), and moreover is a complete section.
We now define a section
\begin{equation*}
    \mathcal{A}^\theta(\mtt{i},\umtt{j})=\{\mtt{k}\in\mathcal{I}^*:b_{\mtt{k}}<b_{\mtt{i}}^{1/\theta-1}\leq b_{\mtt{k}^-}\quad\text{and}\quad\eta(\mtt{k})\preccurlyeq\umtt{j}\}.
\end{equation*}
In words, this set codes the cylinders intersecting $Q$ with height $(\diam Q)^{1/\theta}$ and width greater than $\diam Q$.
However, this set of cylinders may not entirely cover $Q$, so we add the missing cylinders to form
\begin{equation*}
    \mathcal{B}^\theta(\mtt{i},\umtt{j})=\mathcal{A}^\theta(\mtt{i},\umtt{j})\wedge\eta^{-1}(\umtt{j}).
\end{equation*}

Note that if $\mtt{k}\in\mathcal{B}^\theta(\mtt{i},\umtt{j})$, there is a unique $\umtt{l}(\mtt{k})\in\eta(\mathcal{I}^*)$ such that $\umtt{j}=\eta(\mtt{k})\umtt{l}(\mtt{k})$.
Thus the section $\mathcal{B}^\theta(\mtt{i},\umtt{j})$ induces a decomposition of $Q$ into wide pseudo-cylinders:
\begin{equation*}
    Q=P(\mtt{i},\umtt{j})=\bigcup_{\mtt{k}\in\mathcal{B}^\theta(\mtt{i},\umtt{j})}P(\mtt{i}\mtt{k},\umtt{l}(\mtt{k})).
\end{equation*}
By definition of $\mathcal{B}^\theta(\mtt{i},\umtt{j})$, if $\umtt{l}(\mtt{k})\neq\varnothing$, then the height of the corresponding pseudo-cylinder $P(\mtt{i}\mtt{k},\umtt{l}(\mtt{k}))$ is (up to a constant multiple) $b_{\mtt{i}}^{1/\theta}$.

Another equivalent way to think about the decomposition is as follows.
First, observe that
\begin{equation*}
    \mathcal{B}^\theta_{\mathrm{thick}}(\mtt{i},\umtt{j})\coloneqq\mathcal{B}^\theta(\mtt{i},\umtt{j})\cap\eta^{-1}(\umtt{j})=\{\mtt{k}\in\mathcal{B}^\theta(\mtt{i},\umtt{j}):\umtt{l}(\mtt{k})=\varnothing\}.
\end{equation*}
Then set $\mathcal{B}^\theta_{\mathrm{thin}}(\mtt{i},\umtt{j})\coloneqq \mathcal{B}^\theta(\mtt{i},\umtt{j})\setminus\mathcal{B}^\theta_{\mathrm{thick}}(\mtt{i},\umtt{j})$.
The pseudo-cylinders corresponding to the elements of $\mathcal{B}^\theta_{\mathrm{thin}}(\mtt{i},\umtt{j})$ are precisely formed by ``grouping'' the cylinders composing the approximate square $Q$ which have height less than $R^{1/\theta}$ into pseudo-cylinders with height approximately $R^{1/\theta}$.
Moreover, the pseudo-cylinders corresponding to the elements of $\mathcal{B}^\theta_{\mathrm{thick}}(\mtt{i},\umtt{j})$ are in fact cylinders, and they have height greater than $R^{1/\theta}$.

\subsubsection{Defining types and counting type classes}
We first define the notion of the \emph{type} corresponding to a word.
Suppose $\mtt{i}=(i_1,\ldots,i_n)\in\mathcal{I}^n$ for some $n\in\N$, and write
\begin{equation*}
    \bm{\xi}(\mtt{i})=(p_j)_{j\in\mathcal{I}}\qquad\text{where}\qquad p_j=\frac{\#\{k=1,\ldots,n:i_k=j\}}{n}.
\end{equation*}
Of course, $\bm{\xi}(\mtt{i})\in\mathcal{P}$.

Now, fix $\theta\in(0,1)$, an approximate square $Q=P(\mtt{i},\umtt{j})$, and corresponding section $\mathcal{B}^\theta(\mtt{i},\umtt{j})$ as defined in the previous section.
Fix $\mtt{k}\in\mathcal{B}^\theta(\mtt{i},\umtt{j})$ and define the \emph{type} of $\mtt{k}$ as the pair
\begin{equation*}
    \bm{\zeta}(\mtt{k})=(\bm{\xi}(\mtt{i}),\bm{\xi}(\mtt{k}))
\end{equation*}
and denote the set of all types
\begin{equation*}
    \mathcal{T}^\theta(\mtt{i},\umtt{j})=\{\bm{\zeta}(\mtt{k}):\mtt{k}\in\mathcal{B}^\theta(\mtt{i},\umtt{j})\}.
\end{equation*}
Conversely, given a type $(\bm{v},\bm{w})\in\mathcal{T}^\theta(\mtt{i},\umtt{j})$, define the corresponding \emph{type class} by
\begin{equation*}
    \mathcal{C}^\theta(\bm{v},\bm{w})=\mathcal{C}^\theta_{\mtt{i},\umtt{j}}(\bm{v},\bm{w})\coloneqq \left\{\mtt{k}\in\mathcal{B}^{\theta}(\mtt{i},\umtt{j}):\bm{\zeta}(\mtt{k})=(\bm{v},\bm{w})\right\}.
\end{equation*}
An important observation is that if $\mtt{k},\mtt{k}'\in\mathcal{C}^\theta_{\mtt{i},\umtt{j}}(\bm{v},\bm{w})$ have the same type, then
\begin{equation*}
    |\mtt{k}|=|\mtt{k}'|,\quad\eta(\mtt{k})=\eta(\mtt{k}'),\quad a_{\mtt{k}}=a_{\mtt{k}'},\quad\text{and}\quad b_{\mtt{k}}=b_{\mtt{k}'}.
\end{equation*}
We will require the following key estimates on the exponential growth rate of the number of possible types and the size of each type class.
\begin{lemma}\label{l:type-count}
    Fix $\theta\in(0,1)$.
    Then the following hold:
    \begin{enumerate}[nl,r]
        \item\label{im:count} We have
            \begin{equation*}
                \frac{\log\#\mathcal{T}^\theta(\mtt{i},\umtt{j})}{\log(1/b_{\mtt{i}})}=O\left(\frac{\log|\mtt{i}|}{|\mtt{i}|}\right).
            \end{equation*}
        \item\label{im:entropy} Let $\mtt{k}\in\mathcal{B}^\theta(\mtt{i},\umtt{j})$ have type $\bm{\zeta}(\mtt{k})=(\bm{v},\bm{w})$.
            Then
            \begin{equation*}
                \frac{\log\#\mathcal{C}^\theta(\bm{v},\bm{w})}{|\mtt{k}|}=H(\bm{w})-H(\eta(\bm{w}))+O\left(\frac{\log|\mtt{i}|}{|\mtt{i}|}\right).
            \end{equation*}
    \end{enumerate}
\end{lemma}
\begin{proof}
    To see \cref{im:count}, there is a constant $M>0$ so that $|\mtt{k}|\leq M\cdot|\mtt{i}|$ for all $\mtt{k}\in\mathcal{B}^\theta(\mtt{i},\umtt{j})$.
    Thus by \cite[Lemma~2.1.2]{zbl:1177.60035},
    \begin{equation*}
        1\leq \#\mathcal{T}^\theta(\mtt{i},\umtt{j})\leq\sum_{n=0}^{M|\mtt{i}|}\#\{\bm{\xi}(\mtt{k}):\mtt{k}\in\mathcal{I}^{n}\}\leq (M|\mtt{i}|+1)^{\#\mathcal{I}+1}.
    \end{equation*}
    But $\log(1/b_{\mtt{i}})\approx |\mtt{i}|\approx |\mtt{k}|$, from which the result follows.

    Next, we see \cref{im:entropy}.
    Let $\mtt{k}\in\mathcal{B}^\theta(\mtt{i},\umtt{j})$ have type $\bm{\zeta}(\mtt{k})=(\bm{v},\bm{w})$, and recall that for all $\mtt{k}'\in\mathcal{C}^\theta(\bm{v},\bm{w})$,
    \begin{equation*}
        m\coloneqq|\mtt{k}|=|\mtt{k}'|\qquad\text{and}\qquad \eta(\mtt{k})=\eta(\mtt{k}').
    \end{equation*}
    Now, by \cite[Lemma~2.1.8]{zbl:1177.60035},
    \begin{equation*}
        (m+1)^{-\#\mathcal{I}}\exp(m\cdot H(\bm{w}))\leq\#\{\mtt{j}\in\mathcal{I}^m:\bm{\xi}(\mtt{j})=\bm{w}\}\leq \exp(m\cdot H(\bm{w})).
    \end{equation*}
    However, $\mathcal{C}^\theta(\bm{v},\bm{w})$ consists only of those $\mtt{j}$ for which $\eta(\mtt{j})=\eta(\mtt{k})$.
    Moreover, the quantity
    \begin{equation*}
        \#\{\mtt{j}\in\mathcal{I}^m:\bm{\xi}(\mtt{j})=\bm{w}\text{ and }\eta(\mtt{j})=\umtt{h}\}
    \end{equation*}
    is independent of the choice of $\umtt{h}\in\eta(\mathcal{I}^m)$ as long as $\bm{\xi}(\umtt{h})=\eta(\bm{w})$, and 0 otherwise.
    Thus again applying \cite[Lemma~2.1.8]{zbl:1177.60035} to count the number of possible choices for $\umtt{h}$,
    \begin{align*}
        \frac{(m+1)^{-\#\mathcal{I}}\exp(m\cdot H(\bm{w}))}{\exp(m\cdot H(\eta(\bm{w})))}&\leq\#\{\mtt{j}\in\mathcal{I}^m:\bm{\xi}(\mtt{j})=\bm{w}\text{ and }\eta(\mtt{j})=\eta(\mtt{k})\}\\
                                                                               &\leq \frac{\exp(m\cdot H(\bm{w}))}{(m+1)^{-\#\eta(\mathcal{I})}\exp(m\cdot H(\eta(\bm{w})))}.
    \end{align*}
    Taking logarithms, dividing by $m$, and recalling again that $|\mtt{k}|\approx|\mtt{i}|$ yields the desired result.
\end{proof}

\subsection{Covering thin cylinders}\label{ss:thin}
We now begin our covering arguments for the individual types, beginning with the thin cylinders.

We first require a covering lemma for wide pseudo-cylinders in terms of approximate squares.
Recall that $\mathcal{S}$ denotes the set of all approximate squares, and $\mathcal{S}(r)$ denotes the approximate squares with diameter approximately $r$.
If $P(\mtt{i},\umtt{j})$ is a wide pseudo-cylinder, we can write it as a union of the approximate squares in the family
\begin{equation*}
    \mathcal{Q}(\mtt{i},\umtt{j})\coloneqq \{Q\in\mathcal{S}:Q=P(\mtt{i},\umtt{k})\text{ for some }\umtt{k}\in\eta(\mathcal{I}^*)\text{ and }Q\subset P(\mtt{i},\umtt{j})\}.
\end{equation*}
Note that since each $Q=P(\mtt{i},\umtt{k})$ for some $\umtt{k}$, we have $\mathcal{Q}(\mtt{i},\umtt{j})\subset\mathcal{S}(b_{\mtt{i}})$.
We then have the following standard covering result which is given explicitly in, for example, \cite[Lemma~4.8]{arxiv:2309.11971v1}.
We give the short proof for the convenience of the reader.
\begin{lemma}\label{l:pseudo-cyl-cover}
    Let $P(\mtt{i},\umtt{j})$ be a wide pseudo-cylinder.
    Then
    \begin{equation*}
        \#\mathcal{Q}(\mtt{i},\umtt{j})\approx\left(\frac{a_{\mtt{i}\umtt{j}}}{b_{\mtt{i}}}\right)^{\dimB\eta(K)}.
    \end{equation*}
\end{lemma}
\begin{proof}
    Write $\mathcal{Q}(\mtt{i},\umtt{j})=\{Q_1,\ldots,Q_m\}$, and for each $i=1,\ldots,m$, let $\umtt{k}_i$ be such that $Q_i=P(\mtt{i},\umtt{k}_i)$.
    Since $\{\umtt{k}_1,\ldots,\umtt{k}_m\}\succcurlyeq\{\umtt{j}\}$ is a complete section relative to $\umtt{j}$, writing $s=\dimB\eta(K)$ and recalling that $\eta(K)$ is the attractor of a self-similar IFS satisfying the open set condition,
    \begin{equation}
        \sum_{i=1}^m a_{\umtt{k}_i}^s =a_{\umtt{j}}^s.
    \end{equation}
    But $a_{\mtt{i}\umtt{k}_i}\approx b_{\mtt{i}}$ since each $Q_i$ is an approximate square, so the result follows.
\end{proof}
With this lemma, we can now obtain our main covering bound for thin cylinders.
\begin{proposition}\label{p:thin-cover}
    Fix $\theta\in(0,1)$.
    Then for all approximate squares $Q=P(\mtt{i},\umtt{j})$ and type classes $(\bm{v},\bm{w})=\bm{\zeta}(\mtt{k})$ for some $\mtt{k}\in\mathcal{B}^\theta_{\mathrm{thin}}(\mtt{i},\umtt{j})$, the following hold:
    \begin{enumerate}[nl,r]
        \item\label{im:thin-threshold} We have $\phi(\theta,\bm{v})\leq \Gamma(\bm{w})+O\bigl(|\mtt{i}|^{-1}\bigr)$.
        \item\label{im:thin-bounds} Set $E=\bigcup_{\mtt{k}'\in\mathcal{C}^\theta(\bm{v},\bm{w})}P(\mtt{i}\mtt{k}',\mtt{l}(\mtt{k}'))$.
            Then
            \begin{equation*}
                \frac{\log\#\{Q'\in\mathcal{S}(b_{\mtt{i}}^{1/\theta}):Q'\cap E\neq\varnothing\}}{(1/\theta-1)\log(1/b_{\mtt{i}})}=f_{\mathrm{thin}}(\theta,\bm{v},\bm{w})+O\left(\frac{\log|\mtt{i}|}{|\mtt{i}|}\right).
            \end{equation*}
    \end{enumerate}
\end{proposition}
\begin{proof}
    To see \cref{im:thin-threshold}, first observe by the definition of $\mathcal{B}^\theta_{\mathrm{thin}}$ that $b_{\mtt{i}}^{1/\theta-1}\approx b_{\mtt{k}}$, and since $Q$ is an approximate square, $a_{\mtt{i}}a_{\mtt{k}}\gtrsim b_{\mtt{i}}$.
    In particular,
    \begin{equation}\label{e:theta-rel}
        \left(\frac{1}{\theta}-1\right)\cdot |\mtt{i}|\cdot\chi_2(\bm{v}) = |\mtt{k}|\cdot \chi_2(\bm{w})+O(1)
    \end{equation}
    and
    \begin{equation}\label{e:kvalue}
        |\mtt{i}|\cdot\chi_1(\bm{v})+|\mtt{k}|\cdot\chi_1(\bm{w})\leq |\mtt{i}|\cdot\chi_2(\bm{v})+O(1).
    \end{equation}
    Substituting the value of $|\mtt{k}|$ from \cref{e:kvalue} into \cref{e:theta-rel} and dividing through by $|\mtt{i}|$ yields
    \begin{equation*}
        \frac{1}{\Gamma(\bm{v})}+\left(\frac{1}{\theta}-1\right)\cdot\frac{1}{\Gamma(\bm{w})}\leq 1+O\bigl(|\mtt{i}|^{-1}\bigr).
    \end{equation*}
    Since $\Gamma$ takes values in a compact subinterval of $(1,\infty)$, this is a rearrangement of \cref{im:thin-threshold}.

    To see \cref{im:thin-bounds}, first note that for each $\mtt{k}'\in\mathcal{C}^\theta(\bm{v},\bm{w})$,
    \begin{equation*}
        b_{\mtt{i}\mtt{k}'}\approx b_{\mtt{i}}^{1/\theta}\qquad\text{and}\qquad a_{\mtt{i}\mtt{k}'\umtt{l}(\mtt{k}')}=a_{\mtt{i}\umtt{j}}\approx b_{\mtt{i}}
    \end{equation*}
    and moreover by \cref{e:theta-rel} and \cref{l:type-count}~\cref{im:entropy},
    \begin{equation}\label{e:thin-type-count}
        \frac{\log\#\mathcal{C}^\theta(\bm{v},\bm{w})}{(1/\theta-1)\log(1/b_{\mtt{i}})}=
        \frac{\log\#\mathcal{C}^\theta(\bm{v},\bm{w})}{|\mtt{i}|\cdot (1/\theta-1) \chi_2(\bm{v})}=
        \frac{H(\bm{w})-H(\eta(\bm{w}))}{\chi_2(\bm{w})}+O\left(\frac{\log|\mtt{i}|}{|\mtt{i}|}\right).
    \end{equation}
    Thus by \cref{l:pseudo-cyl-cover},
    \begin{align*}
        \#\{Q'\in\mathcal{S}(b_{\mtt{i}}^{1/\theta}):Q'\cap E\neq\varnothing\}&\approx \sum_{k'\in\mathcal{C}^\theta(\bm{v},\bm{w})}\left(\frac{a_{\mtt{i}\mtt{k}'\umtt{l}(\mtt{k'})}}{b_{\mtt{i}\mtt{k}'}}\right)^{\dimB\eta(K)}\\
                                                                              &\approx \#\mathcal{C}^\theta(\bm{v},\bm{w})\cdot b_{\mtt{i}}^{(1-1/\theta)\cdot\dimB\eta(K)}.
    \end{align*}
    Taking logarithms, dividing through by$(1/\theta-1)\log(1/b_{\mtt{i}})$, and applying \cref{e:thin-type-count} completes the proof.
\end{proof}

\subsection{Covering thick cylinders}\label{ss:thick}
We now obtain our main bounds for thick cylinders.
First, we require the following covering lemma for cylinders by approximate squares with height smaller than the height of the cylinder.
This result is a direct application of (the proof of) \cite[Theorem~2.4]{zbl:0757.28011} and the short proof can be found, for example, in the proof of \cite[Lemma~4.9]{arxiv:2309.11971v1}.
Here, we adapt the notation to align with the application in \cref{p:thick-cover} and give a proof for the convenience of the reader.
Recall that $t_{\min}=\dimB K-\dimB\eta(K)$.
\begin{lemma}\label{l:cyl-bound}
    Suppose $\mtt{i},\mtt{k}'\in\mathcal{I}^*$ are such that $b_{\mtt{i}}^{1/\theta-1}\lesssim b_{\mtt{k}'}$ and $a_{\mtt{i}\mtt{k}'}\approx b_{\mtt{i}}$.
    Then
    \begin{align*}
        \#\{Q\in\mathcal{S}(b_{\mtt{i}}^{1/\theta}):Q\subset[\mtt{i}\mtt{k}']\}\approx\left(\frac{1}{b_{\mtt{i}}}\right)^{(1/\theta-1)\cdot\dimB K}\cdot\left(\frac{1}{b_{\mtt{k}'}}\right)^{-t_{\min}}.
    \end{align*}
\end{lemma}
\begin{proof}
    Fix $\mtt{i},\mtt{k}'\in\mathcal{I}^*$.
    Write $R=b_{\mtt{i}}^{1/\theta-1}/b_{\mtt{k}'}$ and enumerate $\mathcal{S}(R)=\{Q_1,\ldots,Q_m\}$.
    Inspecting the proofs of \cite[Lemmas~2.1, 2.2, \& 2.3]{zbl:0757.28011}, we see that
    \begin{equation*}
        m\approx \left(\frac{b_{\mtt{k}'}}{b_{\mtt{i}}^{1/\theta-1}}\right)^{\dimB K}.
    \end{equation*}
    Moreover, for each $i=1,\ldots,m$, we write $Q_i=P(\mtt{j}_i,\umtt{k}_i)$ for some $\mtt{j}_i\in\mathcal{I}^*$ and $\umtt{k}_i\in\eta(\mathcal{I}^*)$, so that
    \begin{equation*}
        \mathcal{Q}(\mtt{i}\mtt{k}'\mtt{j}_i,\umtt{k}_i)\subset \mathcal{S}(b_{\mtt{i}}^{1/\theta})\qquad\text{and}\qquad[\mtt{i}\mtt{k}']=\bigcup_{i=1}^m\bigcup_{Q\in\mathcal{Q}(\mtt{i}\mtt{k}'\mtt{j}_i,\umtt{k}_i)} Q.
    \end{equation*}
    Thus by \cref{l:pseudo-cyl-cover} applied to each wide pseudo-cylinder $P(\mtt{i}\mtt{j}_i,\umtt{k}_i)$, since $Q_i$ is an approximate square so $a_{\mtt{j}_i\umtt{k}_i}\approx b_{\mtt{j}_i}$, and $a_{\mtt{i}\mtt{k}'}\approx b_{\mtt{i}}$ by assumption,
    \begin{align*}
        \#\{Q\in\mathcal{S}(b_{\mtt{i}}^{1/\theta}):Q\subset[\mtt{i}\mtt{k}']\} &=\sum_{i=1}^m\#\mathcal{Q}(\mtt{i}\mtt{k}'\mtt{j}_i,\umtt{k}_i)\\
                                                                        &\approx\sum_{i=1}^m\left(\frac{a_{\mtt{i}\mtt{k}'\mtt{j}_i\umtt{k}_i}}{b_{\mtt{i}\mtt{k}'\mtt{j}_i}}\right)^{\dimB\eta(K)}\\
                                                                        &\approx\left(\frac{1}{b_{\mtt{i}}}\right)^{(1/\theta-1)\cdot\dimB K}\cdot\left(\frac{1}{b_{\mtt{k}'}}\right)^{-t_{\min}}
    \end{align*}
    as claimed.
\end{proof}
With this lemma in hand, we now obtain our main results concerning thick cylinders.
\begin{proposition}\label{p:thick-cover}
    Fix $\theta\in(0,1)$.
    Then for all approximate squares $Q=P(\mtt{i},\umtt{j})$ and type classes $(\bm{v},\bm{w})=\bm{\zeta}(\mtt{k})$ for some $\mtt{k}\in\mathcal{B}^\theta_{\mathrm{thick}}(\mtt{i},\umtt{j})$, the following hold:
    \begin{enumerate}[nl,r]
        \item\label{im:thick-threshold} We have $\phi(\theta,\bm{v})\geq \Gamma(\bm{w})+O\bigl(|\mtt{i}|^{-1}\bigr)$.
        \item\label{im:thick-bounds} Set $E=\bigcup_{\mtt{k}'\in\mathcal{C}^\theta(\bm{v},\bm{w})}P(\mtt{i}\mtt{k}',\mtt{l}(\mtt{k}'))$.
            Then
            \begin{equation*}
                \frac{\log\#\{Q'\in\mathcal{S}(b_{\mtt{i}}^{1/\theta}):Q'\cap E\neq\varnothing\}}{(1/\theta-1)\log(1/b_{\mtt{i}})}=f_{\mathrm{thick}}(\theta,\bm{v},\bm{w})+O\left(\frac{\log|\mtt{i}|}{|\mtt{i}|}\right).
            \end{equation*}
    \end{enumerate}
\end{proposition}
\begin{proof}
    First, by the definition of $\mathcal{B}^\theta_{\mathrm{thick}}(\mtt{i},\umtt{j})$, $b_{\mtt{i}}^{1/\theta-1}\lesssim b_{\mtt{k}}$ and since $Q$ is an approximate square, $a_{\mtt{i}}a_{\mtt{k}}\approx b_{\mtt{i}}$.
    Therefore similar arguments as used in the proof of \cref{p:thin-cover}~\cref{im:thin-threshold} yield \cref{im:thick-threshold}.

    Next, we see \cref{im:thick-bounds}.
    Let $\mtt{k}'\in\mathcal{C}^\theta(\bm{v},\bm{w})$ be arbitrary.
    Since $a_{\mtt{k}'}\approx b_{\mtt{i}}a_{\mtt{i}}^{-1}$, taking logarithms and rearranging gives that
    \begin{equation*}
        |\mtt{i}|\cdot (1/\theta-1)\chi_2(\bm{v})=|\mtt{k}'|\cdot \chi_1(\bm{w})\phi(\theta,\bm{v})+O(1).
    \end{equation*}
    In particular,
    \begin{equation}\label{e:mn-ratio}
        \frac{\log\left( ( 1/b_{\mtt{k}'} )^{-t_{\min}}\right)}{(1/\theta-1)\log(1/b_{\mtt{i}})} = -\frac{t_{\min}}{1/\theta-1}\cdot\frac{|\mtt{k}'|\cdot\chi_2(\bm{w})}{|\mtt{i}|\cdot\chi_2(\bm{v})}
        = -\frac{t_{\min}\chi_2(\bm{w})}{\phi(\theta,\bm{v})\chi_1(\bm{w})}+O(|\mtt{i}|^{-1}),
    \end{equation}
    and by \cref{l:type-count}~\cref{im:entropy},
    \begin{equation}\label{e:type-count}
        \frac{\log\#\mathcal{C}^\theta(\bm{v},\bm{w})}{(1/\theta-1)\log(1/b_{\mtt{i}})}
        =\frac{\log\#\mathcal{C}^\theta(\bm{v},\bm{w})}{(1/\theta-1)|\mtt{i}|\chi_2(\bm{v})}
        =\frac{H(\bm{w})-H(\eta(\bm{w}))}{\phi(\theta,\bm{v})\chi_1(\bm{w})}+O\left(\frac{\log|\mtt{i}|}{|\mtt{i}|}\right).
    \end{equation}
    Moreover, since $\mtt{k}'\in\mathcal{B}^\theta_{\mathrm{thick}}$, we have $\umtt{l}(\mtt{k}')=\varnothing$ so we may apply \cref{l:cyl-bound} to each cylinder $P(\mtt{i}\mtt{k}',\umtt{l}(\mtt{k}'))=[\mtt{i}\mtt{k}']$ giving
    \begin{align*}
        \#\{Q'\in\mathcal{S}(b_{\mtt{i}}^{(1/\theta)}):Q'\cap E\neq\varnothing\}
        &\approx\left(\frac{1}{b_{\mtt{i}}}\right)^{(1/\theta-1)\cdot\dimB K}\cdot\sum_{\mtt{k}'\in\mathcal{C}^\theta(\bm{v},\bm{w})}\left(\frac{1}{b_{\mtt{k}'}}\right)^{-t_{\min}}\\
        &=\left(\frac{1}{b_{\mtt{i}}}\right)^{(1/\theta-1)\cdot\dimB K}\cdot\#\mathcal{C}^\theta(\bm{v},\bm{w})\left(\frac{1}{b_{\mtt{k}}}\right)^{-t_{\min}}.
    \end{align*}
    Recalling the computations in \cref{e:mn-ratio} and \cref{e:type-count}, taking logarithms and dividing by $(1/\theta-1)\log(1/b_\mtt{i})$ gives
    \begin{align*}
        &\frac{\log\#\{Q'\in\mathcal{S}(b_{\mtt{i}}^{(1/\theta)}):Q'\cap E\neq\varnothing\}}{(1/\theta-1)\log(1/b_{\mtt{i}})}\\
        &\qquad = \dimB K+\frac{H(\bm{w})-H(\eta(\bm{w}))}{\phi(\theta,\bm{v})\chi_1(\bm{w})}+\frac{-t_{\min}\chi_2(\bm{w})}{\phi(\theta,\bm{v})\chi_1(\bm{w})}+O\left(\frac{\log|\mtt{i}|}{|\mtt{i}|}\right)\\
        &\qquad =f_{\mathrm{thick}}(\theta,\bm{v},\bm{w})+O\left(\frac{\log|\mtt{i}|}{|\mtt{i}|}\right)
    \end{align*}
    as claimed.
\end{proof}

\subsection{Combining bounds and completing the proof}
To complete the proof, it simply remains to combine the results established in the previous sections.
\begin{proofref}{t:gl-variational}
    Let $\theta\in(0,1)$ be fixed.
    First, there is a constant $M>0$ such that any ball $B(x,r)$ can be covered by $M$ approximate squares in $\mathcal{S}(r)$ and vice versa.
    Thus if we set for $n\in\N$
    \begin{equation*}
        \mathcal{D}_n\coloneqq \left\{P(\mtt{i},\umtt{j}):P(\mtt{i},\umtt{j})\text{ is an approximate square with }|\mtt{i}|=n \right\},
    \end{equation*}
    then
    \begin{equation}\label{e:spec-limsup}
        \dimAs\theta K=\limsup_{n\to\infty}\max_{P(\mtt{i},\umtt{j})\in\mathcal{D}_n}\frac{\log\#\{Q\in\mathcal{S}(b_{\mtt{i}}^{1/\theta}):Q\cap P(\mtt{i},\umtt{j})\neq\varnothing\}}{(1/\theta-1)\log(1/b_{\mtt{i}})}.
    \end{equation}
    Now, fix an approximate square $P(\mtt{i},\umtt{j})$, which we recall that we can decompose as
    \begin{align*}
        P(\mtt{i},\umtt{j})=\bigcup_{\mtt{k}\in\mathcal{B}^{\theta}(\mtt{i},\umtt{j})}P(\mtt{i}\mtt{k},\umtt{l}(\mtt{k})) =\bigcup_{(\bm{v},\bm{w})\in\mathcal{T}^\theta(\mtt{i},\umtt{j})}\bigcup_{\mtt{k}\in\mathcal{C}^\theta_{\mtt{i},\umtt{j}}(\bm{v},\bm{w})}P(\mtt{i}\mtt{k},\umtt{l}(\mtt{k})).
    \end{align*}
    For $(\bm{v},\bm{w})\in\mathcal{T}^\theta(\mtt{i},\umtt{j})$, write
    \begin{equation*}
        N(\bm{v},\bm{w})=\#\left\{Q\in\mathcal{S}(b_{\mtt{i}}^{1/\theta}):Q\cap P(\mtt{i}\mtt{k},\umtt{l}(\mtt{k})) \neq \varnothing \text{ for some }\mtt{k}\in\mathcal{C}^\theta_{\mtt{i},\umtt{j}}(\bm{v},\bm{w})\right\}.
    \end{equation*}

    Next, let $(\bm{v}_0,\bm{w}_0)=\bm{\zeta}(\mtt{k})$ be chosen so that $N(\bm{v}_0,\bm{w}_0)$ is maximised.
    Suppose that $\mtt{k}\in\mathcal{B}_{\mathrm{thin}}^\theta(\mtt{i},\umtt{j})$.
    Then by \cref{l:type-count}~\cref{im:count} and \cref{p:thin-cover}~\cref{im:thin-bounds},
    \begin{align*}
        \frac{\log\#\{Q\in\mathcal{S}(b_{\mtt{i}}^{1/\theta}):Q\cap P(\mtt{i},\umtt{j})\neq\varnothing\}}{(1/\theta-1)\log(1/b_{\mtt{i}})}
        &=\frac{\log N(\bm{v}_0,\bm{w}_0)}{(1/\theta-1)\log(1/b_{\mtt{i}})}+O\left(\frac{\log |\mtt{i}|}{|\mtt{i}|}\right)\\
        &= f_{\mathrm{thin}}(\theta,\bm{v}_0,\bm{w}_0)+O\left(\frac{\log |\mtt{i}|}{|\mtt{i}|}\right).
    \end{align*}
    Moreover, since $f$ and $f_{\mathrm{thin}}$ are continuous functions on the compact domain $\mathcal{P}\times\mathcal{P}$, they are in fact uniformly continuous.
    In particular, for all $\varepsilon>0$ and all $|\mtt{i}|$ sufficiently large depending on $\varepsilon$, by \cref{p:thin-cover}~\cref{im:thin-threshold}, for all $\mtt{k}\in\mathcal{B}_{\mathrm{thin}}^\theta(\mtt{i},\umtt{j})$ with type $(\bm{v},\bm{w})=\bm{\zeta}(\mtt{k})$,
    \begin{equation*}
        |f(\theta,\bm{v},\bm{w})-f_{\mathrm{thin}}(\theta,\bm{v},\bm{w})|\leq\varepsilon.
    \end{equation*}

    Of course, the analogous results hold with $f_{\mathrm{thick}}$ in place of $f_{\mathrm{thin}}$ for $\mtt{k}\in\mathcal{B}_{\mathrm{thick}}^\theta(\mtt{i},\umtt{j})$ by \cref{p:thick-cover}~\cref{im:thick-threshold} and \cref{im:thick-bounds}.

    Combining these observations, for all $\varepsilon>0$ and $n$ sufficiently large depending on $\varepsilon$, with
    \begin{equation*}
        \mathcal{T}^\theta_n\coloneqq \bigcup_{P(\mtt{i},\umtt{j})\in\mathcal{D}_n}\mathcal{T}^\theta(\mtt{i},\umtt{j})
    \end{equation*}
    denoting the set of all types at level $n$,
    \begin{equation*}
        \left\lvert\max_{P(\mtt{i},\umtt{j})\in\mathcal{D}_n}\frac{\log\#\{Q\in\mathcal{S}(b_{\mtt{i}}^{1/\theta}):Q\cap P(\mtt{i},\umtt{j})\neq\varnothing\}}{(1/\theta-1)\log(1/b_{\mtt{i}})}-\max_{(\bm{v},\bm{w})\in\mathcal{T}^\theta_n}f(\theta,\bm{v},\bm{w})\right\rvert\leq\varepsilon.
    \end{equation*}
    But $\mathcal{T}^\theta_n$ becomes arbitrarily dense in $\mathcal{P}\times\mathcal{P}$ as $n$ diverges to infinity, and since $\varepsilon>0$ was arbitrary, the result follows from \cref{e:spec-limsup}.
\end{proofref}

\section{An explicit formula for the Assouad spectrum}\label{s:explicit}
In this section, we obtain the explicit formula for the Assouad spectrum as stated in \cref{it:spectrum-formula}.
Our main tool to solve the optimisation problem which arises in the variational principle in \cref{t:gl-variational} is the duality theory of constrained optimisation.
This approach is based on the general strategy outlined in \cite[§3.1 and §4]{arxiv:2312.08974}, and we recall the main components that we require in \cref{ss:optimisation}.
Since our optimisation problem is only piecewise smooth and far from being (even locally) convex, we also introduce in \cref{ss:projection} some more general conditions which enable the reduction of a certain class of optimisation problem to boundaries of the constraint domain.
With these preliminaries out of the way, in \cref{ss:fibred} we obtain explicit solutions to the optimisation problems that will be relevant for our derivation of the main formula.
The proof of the main result is then completed in \cref{ss:variational}.

\subsection{The geometry of constrained optimisation}\label{ss:optimisation}
In this section, we recall the setup and collect the relevant results from \cite[§3.1]{arxiv:2312.08974}.

We begin by recalling some general definitions from convex analysis.
For a general function $g\colon\R\to\R\cup\{-\infty\}$, the \emph{concave conjugate} is given by
\begin{equation*}
    g^*(\alpha) = \inf_{t\in\R}(t\alpha-g(t)).
\end{equation*}
Note that $g^*$ is always concave, and $g^{**}$ is the concave hull of $g$.
In particular, we always have the inequality
\begin{equation}\label{e:subdiff-rel}
    g^*(\alpha)+g(t)\leq\alpha t,
\end{equation}
and moreover the \emph{subdifferential} $\partial g(t)$ is precisely the set of $\alpha$ for which equality holds in \cref{e:subdiff-rel}.

Suppose moreover that $g$ is a concave function.
We then let $\partial^-g(t)$ (resp.\ $\partial^+g(t)$) denote the left (resp.\ right) derivative of $g$ at $t$, which necessarily exist by concavity of $g$.
Equivalently, $\partial g(t)=[\partial^+g(t),\partial^-g(t)]$.
In particular, $g$ is differentiable at $t$ if and only if $\partial g(t)=\{\alpha\}$, in which case $g'(t)=\alpha$.

Now let $\Delta$ be a compact metric space and let $u,v\colon\Delta\to\R$ be continuous.
We consider the \emph{constrained optimisation}
\begin{equation*}
    F(\alpha)=\max_{\bm{w}\in\Delta}\left\{v(\bm{w}):u(\bm{w})=\alpha\right\},
\end{equation*}
where we set $F(\alpha) = -\infty$ if there does not exist $\bm{w} \in \Delta$ with $u(\bm{w}) = \alpha$.
Associated with this constrained optimisation problem is the \emph{unconstrained dual}
\begin{equation*}
    T(t)=\min_{\bm{w}\in\Delta}\left\{t\cdot u(\bm{w})-v(\bm{w})\right\}.
\end{equation*}
Since $T$ is a minimum of affine functions, $T$ is necessarily concave.
For each $t\in\R$, we denote the set of minimising vectors for $T(t)$ by
\begin{equation*}
    M(t) \coloneqq \left\{\bm{w}\in\Delta:t\cdot u(\bm{w})-v(\bm{w})=T(t)\right\}.
\end{equation*}
The following facts will be useful to us; the proofs of these results are short and elementary and can be found in \cite{arxiv:2312.08974}.
\begin{lemma}[\cite{arxiv:2312.08974}]\label{l:duality}
    The following hold.
    \begin{enumerate}[nl,r]
        \item\label{im:subdiff-char} We have
            \begin{equation*}
                \min_{\bm{w}\in M(t)}u(\bm{w})=\partial^+T(t)\qquad\text{and}\qquad\max_{\bm{w}\in M(t)}u(\bm{w})=\partial^-T(t).
            \end{equation*}
        \item\label{im:conn} Suppose $t\in\R$ and $M(t)$ is connected.
            Then $u(M(t))=\partial T(t)$ and $F(\alpha)=T^*(\alpha)$ for all $\alpha\in\partial T(t)$.
    \end{enumerate}
\end{lemma}
\begin{proof}
    Here, \cref{im:subdiff-char} follows from the proof of \cite[Lemma~3.1]{arxiv:2312.08974} and \cref{im:conn} follows from the proof of \cite[Corollary 3.6]{arxiv:2312.08974}.
\end{proof}
\begin{remark}\label{r:bounded-derivs}
    A straightforward consequence of \cref{l:duality} is that the left and right derivatives of $T$ are elements of the compact set $u(\Delta)\subset\R$.
\end{remark}

\subsection{Attaining the optimisation on the boundary}\label{ss:projection}
We now introduce the concept of an island-free function, and use this to establish some general conditions under which certain constrained optimisation problems are attained on the boundary.
\begin{definition}
    Let $\Delta$ be a topological space and let $f\colon\Delta\to\R$.
    We say that $f$ is \emph{island-free} if for all $t\in\R$ the upper level set
    \begin{equation*}
        \{x\in\Delta:f(x)\geq t\}
    \end{equation*}
    is connected.
\end{definition}
Here, the empty set is always connected.
The following lemma provides some simple conditions to establish island-freeness.
\begin{lemma}\label{l:qc-const}
    The following hold.
    \begin{enumerate}[nl,r]
        \item\label{im:f-quot} Let $\Delta$ be a convex space.
            Suppose $f\colon\Delta\to\R$ is concave and $g\colon\Delta\to(0,\infty)$ is affine.
            Then $f/g$ is island-free.
        \item\label{im:f-prod} Suppose $\Delta_i$ is a topological space and $f_i\colon\Delta_i\to\R_{\geq 0}$ is island-free for each $i=1,\ldots,m$.
            Equip $\Delta\coloneqq\Delta_1\times\cdots\times\Delta_m$ with the product topology.
            Then $f\colon\Delta\to\R_{\geq 0}$ defined by
            \begin{equation*}
                f(x_1,\ldots,x_m)=f_1(x_1)\cdots f_m(x_m)
            \end{equation*}
            is island-free.
    \end{enumerate}
\end{lemma}
\begin{proof}
    To see \cref{im:f-quot}, let $f$ be concave and $g$ affine.
    Then for each $t\in\R$, using positivity of $g$,
    \begin{equation*}
        \{x\in\Delta:f(x)/g(x)\geq t\}=\{x\in\Delta:f(x)-t g(x)\geq 0\}.
    \end{equation*}
    This is a convex set (and in particular connected) since $f(x)-t g(x)$ is concave as a function of $x$.

    Next, we see \cref{im:f-prod}.
    Let $\Delta_i$ and $f_i$ be defined as in the statement of the lemma and let $t\in\R$ be fixed.
    Set
    \begin{equation*}
        A(t)\coloneqq\left\{(x_1,\ldots,x_m):f_1(x_1)\cdots f_m(x_m)\geq t\right\}.
    \end{equation*}
    We must show that $A(t)$ is connected.
    First, let $(y_1,\ldots,y_m)\in A(t)$ be arbitrary.
    Then
    \begin{equation*}
        E(y_1,\ldots,y_m)\coloneqq \prod_{i=1}^m\{x_i:f_i(x_i)\geq f_i(y_i)\}\subset A(t)
    \end{equation*}
    is a connected set since each $f_i$ is island-free.
    Moreover, if $(z_1,\ldots,z_m)\in A(t)$ is arbitrary, since the $f_i$ are non-negative,
    \begin{equation*}
        \varnothing\neq\prod_{i=1}^m\bigl\{x_i:f_i(x_i)\geq\max\{f_i(y_i),f_i(z_i)\}\bigr\}\subset E(y_1,\ldots,y_m)\cap E(z_1,\ldots,z_m).
    \end{equation*}
    Since the union of two intersecting connected sets is connected, any two elements of $A(t)$ are contained in a connected subset of $A(t)$.
    Thus $A(t)$ is connected.
\end{proof}
The preceding lemma, along with concavity of $\bm{w}\mapsto H(\bm{w})-H(\eta(\bm{w}))$ (which follows by the log-sum inequality, see for instance \cite[§2.7]{zbl:1140.94001}), yields the following result.
\begin{corollary}\label{c:f-qc}
    For all $\theta\in(0,1)$, the functions $(\bm{v},\bm{w})\mapsto f_{\mathrm{thin}}(\theta,\bm{v},\bm{w})$ and $(\bm{v},\bm{w})\mapsto f_{\mathrm{thick}}(\theta,\bm{v},\bm{w})$ are island-free.
\end{corollary}
Our main use for island-freeness appears in the following elementary lemma, which provides a partial description of the constrained maximisers of an island-free function in the case that the constrained maximum is not equal to the global maximum.
\begin{lemma}\label{l:qc-proj}
    Let $\Delta$ be a compact metric space and let $f\colon\Delta\to\R$ be continuous and island-free.
    Then the set of global maximisers
    \begin{equation*}
        M\coloneqq\left\{x\in\Delta:f(x)=\max_{y\in\Delta}f(y)\right\}
    \end{equation*}
    is a non-empty, compact and connected subset of $\Delta$.

    Moreover, suppose $E\subset\Delta$ is a non-empty compact set satisfying $M\cap(\Delta\setminus E)\neq\varnothing$.
    Then the set of constrained maximisers
    \begin{equation*}
        M_E\coloneqq\left\{x\in E:f(x)=\max_{y\in E}f(y)\right\},
    \end{equation*}
    intersects the topological boundary of $E$.
\end{lemma}
\begin{proof}
	Since $f$ is continuous, it is immediate that $M$ is well-defined, non-empty and compact.
	Since $f$ is island-free, $M$ is connected.

    Next, suppose $E\subset\Delta$ is compact and $M\cap(\Delta\setminus E)\neq\varnothing$.
    Consider the set
    \begin{equation*}
        G\coloneqq \left\{x\in\Delta:f(x)\geq \max_{y\in E}f(y)\right\}.
    \end{equation*}
    Note that $M_E=G\cap E$ so $G\cap E\neq\varnothing$, and $M\subset G$ so $G\cap(\Delta\setminus E)\neq\varnothing$.
    But $f$ is island-free, so $G$ is connected and therefore $M_E$ intersects the boundary of $E$ relative to $\Delta$.
\end{proof}

\subsection{Fibred optimisers}\label{ss:fibred}
We now solve a few useful global minimisation problems, and also see how these minimisation problems encode the Assouad dimension of the Gatzouras--Lalley carpet $K$.

For $t\in\R$, write
\begin{equation}\label{e:gt-def}
    g(t)=\min_{\bm{w}\in\mathcal{P}}\left\{\frac{t\chi_2(\bm{w})-[H(\bm{w})-H(\eta(\bm{w}))]}{\chi_1(\bm{w})}\right\}.
\end{equation}
We first prove in \cref{p:minimisers} that this definition of $g(t)$ coincides with the alternative definition given in \cref{e:gt-def-intro}.
Moreover, we will give an explicit description of the set of minimisers.

For each $\jh\in\eta(\mathcal{I})$, recall the definition of $\psi_{\jh}$ from the introduction:
\begin{equation*}
    \psi_{\jh}(t)=\frac{\log \sum_{i\in\eta^{-1}(\jh)} b_i^t}{\log a_{\jh}}.
\end{equation*}
Equivalently, for each $i\in\eta^{-1}(\jh)$,
\begin{equation}\label{e:psi-equiv}
    \frac{b_i^t}{\sum_{\ell\in\eta^{-1}(\jh)}b_\ell^t}= a_{\jh}^{-\psi_{\jh}(t)} b_i^t.
\end{equation}
Therefore given $\bm{p}\in\eta(\mathcal{P})$, we may define a probability vector
\begin{equation*}
    \bm{z}(t,\bm{p}) \coloneqq \left(p_{\eta(i)}\cdot a_{\eta(i)}^{-\psi_{\eta(i)}(t)}b_i^t\right)_{i\in\mathcal{I}}.
\end{equation*}
The reason for introducing $\bm{z}(t,\bm{p})$ will become clear below.

We also recall the definition of the \defn{Kullback--Leibler divergence} of two probability vectors $\bm{w}$ and $\bm{v}$ as
\begin{equation*}
    \DKL{\bm{w}}{\bm{v}}=\sum_{i\in\mathcal{I}}w_i\log\left(\frac{w_i}{v_i}\right),
\end{equation*}
where we set $0\log(0/v_i)=0$ regardless of the value of $v_i$.
In general, $\DKL{\bm{w}}{\bm{v}}\geq 0$ with equality if and only if $\bm{w}=\bm{v}$.

Finally, we introduce some notation to denote the set of minimisers for $g(t)$.
For each $t\in\R$, let
\begin{equation*}
    \mathcal{J}(t)=\Bigl\{\jh\in\eta(\mathcal{I}):\psi_{\jh}(t)=\min_{\ih\in\eta(\mathcal{I})}\psi_{\ih}(t)\Bigr\}.
\end{equation*}
We then write
\begin{equation*}
    \mathcal{R}(t)=\{\bm{p}\in\eta(\mathcal{P}):\supp\bm{p}\subset\mathcal{J}(t)\}
    \qquad\text{and}\qquad
    \mathcal{Z}(t)=\{\bm{z}(t,\bm{p}):\bm{p}\in\mathcal{R}(t)\}.
\end{equation*}
We now have the following formula for the function $g(t)$.
\begin{proposition}\label{p:minimisers}
    For each $t\in\R$, using the definition of $g(t)$ from \cref{e:gt-def}, we have
    \begin{equation*}
        g(t)=\min\{\psi_{\jh}(t): \jh\in\eta(\mathcal{I})\}.
    \end{equation*}
    Moreover, the set of minimising vectors satisfies
    \begin{equation*}
        \left\{\bm{w}\in\mathcal{P}:\frac{t\chi_2(\bm{w})-[H(\bm{w})-H(\eta(\bm{w}))]}{\chi_1(\bm{w})}=g(t)\right\}=\mathcal{Z}(t).
    \end{equation*}
\end{proposition}
\begin{proof}
    Suppose $\bm{w}\in\mathcal{P}$ is arbitrary and let $\eta(\bm{w})=\bm{p}$.
    We then compute
    \begin{align*}
        0 &\leq \DKL{\bm{w}}{\bm{z}(t,\bm{p})}\\
          &= \sum_{\jh\in\eta(\mathcal{I})}\sum_{i\in\eta^{-1}(\jh)} w_i \log\left( w_i p_{\jh}^{-1}a_{\jh}^{\psi_{\jh}(t)} b_i^{-t}\right)\\
          &= -H(\bm{w})+H(\eta(\bm{w}))+\sum_{\jh\in\eta(\mathcal{I})}p_{\jh}\psi_{\jh}(t)\log a_{\jh}+t\chi_2(\bm{w}).
    \end{align*}
    Rearranging, we obtain that
    \begin{align*}
        \frac{t\chi_2(\bm{w})-[H(\bm{w})-H(\eta(\bm{w}))]}{\chi_1(\bm{w})} &\geq \frac{-\sum_{\jh\in\eta(\mathcal{I})}p_{\jh}\psi_{\jh}(t)\log a_{\jh}}{\chi_1(\bm{w})}\\
                                                                           &\geq \min\{\psi_{\jh}(t):\jh\in\eta(\mathcal{I})\}.
    \end{align*}
    Observe that the second equality holds if and only if $\bm{p}\in\mathcal{R}(t)$.
    Moreover, for all $\bm{p}\in\eta(\mathcal{P})$, the first equality holds if and only if $\bm{w}=\bm{z}(t,\bm{p})$.
    Thus the desired result follows.
\end{proof}
To conclude this section, we show how the function $g(t)$ encodes the Assouad dimension.
Recall from the introduction that for $\jh\in\eta(\mathcal{I})$, $s_{\jh}$ is the unique solution to the equation
\begin{equation*}
    \sum_{i\in\eta^{-1}(\jh)}b_i^{s_{\jh}}=1.
\end{equation*}
Moreover, we recall the definition of $t_{\max}$, which by the main result of \cite{zbl:1278.37032} satisfies
\begin{equation*}
    t_{\max}=\max_{\jh\in\eta(\mathcal{I})}s_{\jh}=\dimA K-\dimB\eta(K).
\end{equation*}
It turns out that $t_{\max}$ is precisely the unique zero of $g$.
\begin{lemma}\label{l:fibred-max}
    We have
    \begin{align*}
        \max_{\bm{w}\in\mathcal{P}}\left\{\frac{H(\bm{w})-H(\eta(\bm{w}))}{\chi_2(\bm{w})}\right\}=t_{\max}.
    \end{align*}
    Moreover, $t_{\max}$ is the unique zero of $g$ and the set of maximising probability vectors is $\mathcal{Z}(t_{\max})$.
\end{lemma}
\begin{proof}
    Write
    \begin{equation*}
        \mathcal{R}_{\max}=\Bigl\{\bm{p}\in\eta(\mathcal{P}):\supp\bm{p}\subset\{\jh\in\eta(\mathcal{I}):s_{\jh}=t_{\max}\}\Bigr\}.
    \end{equation*}
    Now suppose $\bm{w}\in\mathcal{P}$ is arbitrary.
    Write $\bm{p}=\eta(\bm{w})$ and $\bm{v}=\bigl(p_{\eta(i)}b_i^{s_{\eta(i)}}\bigr)_{i\in\mathcal{I}}$.
    Then
    \begin{align*}
        0 &\leq \DKL{\bm{w}}{\bm{v}}\\
          &=\sum_{\jh\in\eta(\mathcal{I})}\sum_{i\in\eta^{-1}(\jh)}w_i\log (w_i p_{\jh}^{-1}b_i^{-s_{\jh}})\\
          &=-H(\bm{w})+H(\eta(\bm{w}))-\sum_{\jh\in\eta(\mathcal{I})}s_{\jh}\sum_{i\in\eta^{-1}(\jh)}w_i\log b_i\\
          &\leq -H(\bm{w})+H(\eta(\bm{w}))+t_{\max}\chi_2(\bm{w}).
    \end{align*}
    The second inequality is an equality if and only if $\supp\bm{p}\subset\mathcal{R}_{\max}$, in which case the first inequality is an equality if and only if $\bm{w}=\bm{v}$.
    But if $\bm{w}=\bm{v}$ (and $\supp\bm{p}\subset\mathcal{R}_{\max}$), then $\bm{w}=\bm{z}(t_{\max},\bm{p})$ by the definition of $\bm{v}$.

    We have shown that
    \begin{equation*}
        \frac{t_{\max}\chi_2(\bm{w})-[H(\bm{w})-H(\eta(\bm{w}))]}{\chi_1(\bm{w})}\geq 0
    \end{equation*}
    with equality if and only if $\bm{w}=\bm{z}(t_{\max},\bm{p})$ for some $\bm{p}\in\mathcal{R}_{\max}$.
    But this is precisely the same minimisation as the definition of $g(t_{\max})$, so the set of maximising probability vectors must be $\mathcal{Z}(t_{\max})$ by \cref{p:minimisers}.
    (Alternatively, observe that $\mathcal{R}_{\max}=\mathcal{R}(t_{\max})$.)

    To see that $t_{\max}$ is unique with this property, since the left and right derivatives of $g(t)$ lie in the compact interval $\Gamma(\mathcal{P})\subset(1,\infty)$ (see \cref{r:bounded-derivs}), it follows that $g(t)$ is strictly increasing.
\end{proof}

\subsection{Solving the variational formula}\label{ss:variational}
Finally, we can establish an explicit formula for the Assouad spectrum of $K$ by solving the maximisation process underlying the variational formula.
This proof can be subdivided effectively into three parts:
\begin{enumerate}[nl]
    \item First, solve the unconstrained maximisation problems corresponding to the functions $f_{\mathrm{thick}}$ and $f_{\mathrm{thin}}$ and determine the values of $\theta$ for which the respective unconstrained and constrained maxima agree (this is \cref{l:global}).
    \item Next, recalling that $f_{\mathrm{thick}}=f_{\mathrm{thin}}$ on $\Delta_{\mathrm{thick}}(\theta)\cap\Delta_{\mathrm{thin}}(\theta)$, solve the corresponding boundary maximisation (this is \cref{l:boundary}).
    \item Finally, using island-freeness and \cref{l:qc-proj}, reduce the general maximisation to the above cases.
\end{enumerate}
With this outline in mind, we begin the proof.

Recall that we defined
\begin{equation*}
    \theta_{\min}=\phi^{-1}\bigl(\partial^+ g(t_{\min})\bigr)\qquad\text{and}\qquad\theta_{\max}=\phi^{-1}\bigl(\partial^- g(t_{\max})\bigr),
\end{equation*}
and recall the definitions of $\Delta_{\mathrm{thin}}(\theta)$ and $\Delta_{\mathrm{thick}}(\theta)$ from \cref{e:definedeltathin}.
We begin by applying the results from \cref{ss:fibred} to solve the unconstrained maximisation problems corresponding to $f_{\mathrm{thick}}$ and $f_{\mathrm{thin}}$.
\begin{lemma}\label{l:global}
    Let $\theta\in(0,1)$ be arbitrary.
    Given $t\in\R$, let $P_t\colon\R\to\R\cup\{-\infty\}$ denote the function defined by $P_t(t)=g(t)$ and $P_t(x)=-\infty$ for $x\neq t$.
    \begin{enumerate}[nl,r]
        \item\label{im:thin} We have
            \begin{equation*}
                \max_{(\bm{v},\bm{w})\in\mathcal{P}\times\mathcal{P}}f_{\mathrm{thin}}(\theta,\bm{v},\bm{w})=\dimA K=\dimB\eta(K)+\frac{P_{t_{\max}}^*(\phi(\theta))}{\phi(\theta)}
            \end{equation*}
            and the set of probability vectors for which the maximum is attained is given by
            \begin{equation*}
                E_{\mathrm{thin}}\coloneqq \mathcal{P}\times\mathcal{Z}(t_{\max}).
            \end{equation*}
            In particular, $E_{\mathrm{thin}}\cap\Delta_{\mathrm{thin}}(\theta)\neq\varnothing$ if and only if $\theta\geq\theta_{\max}$.
        \item \label{im:thick} We have
            \begin{equation*}
                \max_{(\bm{v},\bm{w})\in\mathcal{P}\times\mathcal{P}}f_{\mathrm{thick}}(\theta,\bm{v},\bm{w})=\dimB K-\frac{g(t_{\min})}{\phi(\theta)}=\dimB\eta(K)+\frac{P_{t_{\min}}^*(\phi(\theta))}{\phi(\theta)},
            \end{equation*}
            and the set of probability vectors for which the maximum is attained is given by
            \begin{equation*}
                E_{\mathrm{thick}}\coloneqq \Gamma^{-1}(\kappa_{\max})\times\mathcal{Z}(t_{\min}).
            \end{equation*}
            In particular, $E_{\mathrm{thick}}\cap\Delta_{\mathrm{thick}}(\theta)\neq\varnothing$ if and only if $\theta\leq\theta_{\min}$.
    \end{enumerate}
\end{lemma}
\begin{proof}
    To see \cref{im:thin}, by inspecting the definition of $f_{\mathrm{thin}}$, the first equality and the formula for $E_{\mathrm{thin}}$ are immediate consequences of \cref{l:fibred-max}, recalling that $\dimA K=t_{\max}+\dimB\eta(K)$.
    To see the second equality, since $\phi(\theta)\in\partial P_{t_{\max}}(t_{\max})$ and since $P_{t_{\max}}(t_{\max})=0$ by \cref{l:fibred-max},
    \begin{equation*}
        \dimB\eta(K)+\frac{P_{t_{\max}}^*(\phi(\theta))}{\phi(\theta)}=\dimB\eta(K)+\frac{t_{\max}\cdot\phi(\theta)-P_{t_{\max}}(t_{\max})}{\phi(\theta)}=\dimA K
    \end{equation*}
    by the definition of $t_{\max}$.
    Finally, recalling that $\theta_{\max}$ is defined from \cref{e:thetaminmaxdef} and using \cref{p:minimisers} and \cref{l:duality}~\cref{im:subdiff-char},
    \begin{equation*}
        \max_{\bm{w}\in\mathcal{Z}(t_{\max})}\Gamma(\bm{w}) = \partial^- g(t_{\max})  =\phi(\theta_{\max}).
    \end{equation*}
    Recall that $\phi(\theta)=\inf_{\bm{v}\in\mathcal{P}}\phi(\theta,\bm{v})$, and $\phi$ is decreasing in $\theta$.
    Therefore $E_{\mathrm{thin}}\cap\Delta_{\mathrm{thin}}(\theta) \neq\varnothing$ if and only if $\theta\geq\theta_{\max}$.

    To see \cref{im:thick}, the maximisation in $\bm{v}$ is clearly attained by any $\bm{v}$ which is supported on indices for which the logarithmic eccentricity $\Gamma$ is as large as possible.
    The first equality and the formula for $E_{\mathrm{thick}}$ then follows from \cref{p:minimisers} since the remaining term in $f_{\mathrm{thick}}(\theta,\bm{v},\bm{w})$ is precisely the negative of the reciprocal of the objective function defining $g(t_{\min})$.
    To see the second inequality, since $\phi(\theta)\in\partial P_{t_{\min}}(t_{\min})$,
    \begin{equation*}
        P_{t_{\min}}(t_{\min})+P_{t_{\min}}^*(\phi(\theta))=t_{\min}\cdot\phi(\theta).
    \end{equation*}
    But recall that $t_{\min}=\dimB K-\dimB\eta(K)$, so rearranging gives
    \begin{equation*}
        \dimB\eta(K)+\frac{P_{t_{\min}}^*(\phi(\theta))}{\phi(\theta)}=\dimB K-\frac{g(t_{\min})}{\phi(\theta)}
    \end{equation*}
    as claimed.
    Finally, let $(\bm{v},\bm{w})\in\Gamma^{-1}(\kappa_{\max})\times\mathcal{Z}(t_{\min})$ be arbitrary.
    By \cref{p:minimisers} and \cref{l:duality}~\cref{im:subdiff-char},
    \begin{equation*}
        \min_{\bm{w}\in\mathcal{Z}(t_{\min})}\Gamma(\bm{w})=\partial^+g(t_{\min})=\phi(\theta_{\min}).
    \end{equation*}
    Now, $\phi(\theta,\bm{v})=\phi(\theta)$, so $E_{\mathrm{thick}}\cap\Delta_{\mathrm{thick}}(\theta)\neq\varnothing$ if and only if $\theta\leq \theta_{\min}$.
\end{proof}
Next, we solve the maximisation problem constrained to the boundary.
For notational simplicity, we write
\begin{equation*}
    \Delta(\theta)=\Delta_{\mathrm{thin}}(\theta)\cap\Delta_{\mathrm{thick}}(\theta).
\end{equation*}
This is the topological boundary of $\Delta_{\mathrm{thin}}(\theta)$ and $\Delta_{\mathrm{thick}}(\theta)$ (relative to $\mathcal{P}\times\mathcal{P}$).
We also recall that $f=f_{\mathrm{thin}}=f_{\mathrm{thick}}$ on $\Delta(\theta)$.
\begin{lemma}\label{l:boundary}
    Suppose $(\bm{v},\bm{w})\in\Delta(\theta)$.
    If $\theta<\theta_{\max}$, then
    \begin{align*}
        \max_{(\bm{v},\bm{w})\in\Delta(\theta)}f(\theta,\bm{v},\bm{w})
        &=\max_{(\bm{v},\bm{w})\in\Delta_{\mathrm{thin}}(\theta)}f_{\mathrm{thin}}(\theta,\bm{v},\bm{w})\\
        &=\dimB\eta(K)+\frac{g^*(\phi(\theta))}{\phi(\theta)}.
    \end{align*}
\end{lemma}
\begin{proof}
    In order to prove the desired formula on the boundary, first introduce the auxiliary function
    \begin{equation*}
        F(\alpha)=\max_{\bm{w}\in\mathcal{P}}\left\{\frac{H(\bm{w})-H(\eta(\bm{w}))}{\chi_1(\bm{w})}:\Gamma(\bm{w})=\alpha\right\}.
    \end{equation*}
    Note that $F(\alpha)$ is the constrained optimisation problem corresponding to the unconstrained problem $g(t)$.
    Moreover, by \cref{p:minimisers}, the minimisation defining $g(t)$ is attained precisely on the set $\mathcal{Z}(t)$, which is a connected set.
    (Alternatively, connectedness of the set of minimisers can be indirectly observed since the negative of the objective function defining $g(t)$ is island-free by \cref{l:qc-const}.)
    Thus applying \cref{l:duality}~\cref{im:conn}, $F(\alpha)=g^*(\alpha)$ for all $\alpha\in \Gamma(\mathcal{Z}(t))$.

    We now obtain the desired bounds.
    First,
    \begin{equation}\label{e:boundaryorhalfspacemax}
        \begin{split}
            \max_{(\bm{v},\bm{w})\in\Delta(\theta)}&f(\theta,\bm{v},\bm{w})-\dimB\eta(K)\\
                                                   &\leq \max_{(\bm{v},\bm{w})\in\Delta_{\mathrm{thin}}(\theta)}f_{\mathrm{thin}}(\theta,\bm{v},\bm{w})-\dimB\eta(K)\\
                                                   &=\max_{(\bm{v},\bm{w})\in\mathcal{P}\times\mathcal{P}}\left\{\frac{H(\bm{w})-H(\eta(\bm{w}))}{\chi_2(\bm{w})}:\Gamma(\bm{w})\geq\phi(\theta,\bm{v})\right\}\\
                                                   &\leq\max_{\bm{w}\in\mathcal{P}}\left\{\frac{H(\bm{w})-H(\eta(\bm{w}))}{\chi_2(\bm{w})}:\Gamma(\bm{w})\geq\phi(\theta)\right\}
        \end{split}
    \end{equation}
    since $\inf_{\bm{v}\in\mathcal{P}}\phi(\theta,\bm{v})=\phi(\theta)$.
    Now recalling \cref{l:global}~\cref{im:thin}, the unconstrained optimisation
    \begin{equation*}
        \max_{\bm{w}\in\mathcal{P}}\left\{\frac{H(\bm{w})-H(\eta(\bm{w}))}{\chi_2(\bm{w})}\right\}
    \end{equation*}
    is attained precisely on the set $\mathcal{Z}(t_{\max})$ and, since $\theta<\theta_{\max}$, $\Gamma(\bm{w})<\phi(\theta)$ for all $\bm{w}\in\mathcal{Z}(t_{\max})$.
    Thus by \cref{l:qc-proj} (island-freeness again follows by \cref{l:qc-const}), the maximisation is attained on the boundary $\{\bm{w}\in\mathcal{P}:\Gamma(\bm{w})=\phi(\theta)\}$.
    But
    \begin{equation*}
        \Gamma^{-1}(\kappa_{\max})\times\{\bm{w}\in\mathcal{P}:\Gamma(\bm{w})=\phi(\theta)\}\subset\Delta(\theta),
    \end{equation*}
    so in fact there is equality throughout \cref{e:boundaryorhalfspacemax} and
    \begin{align*}
        \max_{(\bm{v},\bm{w})\in\Delta(\theta)}&f(\theta,\bm{v},\bm{w})-\dimB\eta(K)\\
                                               &=\max_{\bm{w}\in\mathcal{P}}\left\{\frac{H(\bm{w})-H(\eta(\bm{w}))}{\chi_2(\bm{w})}:\Gamma(\bm{w})=\phi(\theta)\right\}\\
                                               &=\max_{\bm{w}\in\mathcal{P}}\left\{\frac{1}{\phi(\theta)}\cdot \frac{H(\bm{w})-H(\eta(\bm{w}))}{\chi_1(\bm{w})}:\Gamma(\bm{w})=\phi(\theta)\right\}\\
                                               &=\frac{F(\phi(\theta))}{\phi(\theta)}.
    \end{align*}
    In the second equality, we used the substitution $\chi_2(\bm{w})=\Gamma(\bm{w})\chi_1(\bm{w})$.
    Recalling that $F(\phi(\theta))=g^*(\phi(\theta))$ yields the claimed formula.
\end{proof}
We can finally complete the proof of our main result.
\begin{restatement}{it:spectrum-formula}
    Let $K$ be a Gatzouras--Lalley carpet and let
    \begin{equation*}
        \tau(t)=
        \begin{cases}
            g(t) &: t\in [t_{\min},t_{\max}]\\
            -\infty &:\text{ otherwise}.
        \end{cases}
    \end{equation*}
    Then for all $\theta\in(0,1)$,
    \begin{equation*}
        \dimAs\theta K=\dimB\eta(K)+\frac{\tau^*(\phi(\theta))}{\phi(\theta)}.
    \end{equation*}
\end{restatement}
\begin{proof}
    Let $0<\theta<1$.
    Recall that $f_{\mathrm{thin}}$ and $f_{\mathrm{thick}}$ are island-free by \cref{c:f-qc}.
    Moreover, recall from \cref{t:gl-variational} that
    \begin{equation*}
        \dimAs\theta K = \max_{(\bm{v},\bm{w})\in\mathcal{P}\times\mathcal{P}}f(\theta,\bm{v},\bm{w})
    \end{equation*}
    where $f=f_{\mathrm{thin}}$ on $\Delta_{\mathrm{thin}}(\theta)$ and $f=f_{\mathrm{thick}}$ on $\Delta_{\mathrm{thick}}(\theta)$.

    If $\theta\leq\theta_{\min}$, then by \cref{l:global}, $f_{\mathrm{thin}}$ does not attain its global maximum but $f_{\mathrm{thick}}$ does.
    Therefore by \cref{l:qc-proj} as well as the formula in \cref{l:global}~\cref{im:thick},
    \begin{equation*}
        \max_{(\bm{v},\bm{w})\in\Delta_{\mathrm{thick}}(\theta)}f_{\mathrm{thick}}(\theta,\bm{v},\bm{w})=\dimB\eta(K)+\frac{\tau^*(\phi(\theta))}{\phi(\theta)}\geq\max_{(\bm{v},\bm{w})\in\Delta_{\mathrm{thin}}(\theta)}f_{\mathrm{thin}}(\theta,\bm{v},\bm{w})
    \end{equation*}
    yielding the desired formula.
    The analogous argument provides the result for $\theta\geq\theta_{\max}$.

    Otherwise, suppose $\theta_{\min}<\theta<\theta_{\max}$.
    By \cref{l:global}, the unconstrained maxima are not attained for either $f_{\mathrm{thin}}$ or $f_{\mathrm{thick}}$.
    Therefore by \cref{l:qc-proj},
    \begin{equation*}
        \dimAs\theta K = \max_{(\bm{v},\bm{w})\in\Delta(\theta)}f(\theta,\bm{v},\bm{w})= \dimB\eta(K)+\frac{\tau^*(\phi(\theta))}{\phi(\theta)},
    \end{equation*}
    where in the last equality we applied \cref{l:boundary} and used the fact that $\theta_{\min}<\theta<\theta_{\max}$ so $g^*(\phi(\theta))=\tau^*(\phi(\theta))$.
\end{proof}
\begin{remark}\label{r:optim}
    The choice of $\bm{v}_0$ to maximise $f(\theta,\bm{v}_0,\bm{w})$ does not depend on the choice of $\theta$; we can simply take $\bm{v}_0$ to be any probability vector fully supported on the indices $i$ for which $(\log b_i)/(\log a_i)=\kappa_{\max}$.
    This is the reason for the appearance of $\kappa_{\max}$ in the parameter change $\phi(\theta)$.

    The dependence of $\bm{w}$ on $\theta$ is more complex.
    If $\theta$ is such that $\phi(\theta) = \partial^+ g(t(\theta))$ (resp.\ $\partial^- g(\theta)$) for some value $t(\theta) \in (t_{\min},t_{\max})$, then $\eta(\bm{z}(\theta))$ can be taken to be supported on a single column.
    In this case, by the formula for the optimisation given in \cref{p:minimisers}, $\bm{z}(\theta)$ can be given more explicitly as $\bm{z}(\theta)=\bm{z}(t(\theta),\bm{\delta}_{\jh})$, where $\bm{\delta}_{\jh}$ is the probability vector fully supported on a column satisfying $g(t) = \psi_{\jh}(t)$ for some $\varepsilon>0$ and $t\in(t(\theta), t(\theta)+\varepsilon)$ (resp.\ $(t(\theta)-\varepsilon, t(\theta))$).
    Otherwise, if $\theta<\theta_{\min}$ or $\theta>\theta_{\max}$, then by \cref{l:global} we can again take $\bm{z}(\theta)$ to be an explicit vector supported on a single column.
    Finally, if $\partial g^-(t(\theta)) < \phi(\theta) < \partial g^+ (t(\theta))$ for some $t(\theta)$, then $\bm{z}(\theta)$ can be taken to be a convex combination of the optimising vectors corresponding to $\partial g^-(t(\theta))$ and $\partial g^+(t(\theta))$, so $\eta(\bm{z}(\theta))$ can be taken to be supported on at most two columns.
\end{remark}

\begin{acknowledgements}
    AB, JMF and IK were supported by a Leverhulme Trust Research Project Grant (RPG-2019-034) at the University of St Andrews.
    AB was also supported by an EPSRC New Investigators Award (EP/W003880/1) while at Loughborough University.
    AR was supported by EPSRC Grant EP/V520123/1 and the Natural Sciences and Engineering Research Council of Canada.
    IK was also supported by the European Research Council Marie Skłodowska-Curie Actions Postdoctoral Fellowship \#101109013.
    The authors thank Balázs Bárány and Lars Olsen for pointing out a mistake in an earlier version of the proof of \cref{p:minimisers}.
\end{acknowledgements}
\end{document}

%% file: figures/carpet_maps/fig.tex
\begin{tikzpicture}[scale=5]
    \begin{scope}[font=\tiny]
        \node[below] at (0,0) {$0$};
        \node[below] at (1,0) {$1$};
        \node[left] at (0,0) {$0$};
        \node[left] at (0,1) {$1$};
    \end{scope}
    \draw (0,0) rectangle (1,1);

    \draw[dashed] (1/3,0) -- (1/3,1);
    \draw[dashed] (1/2,0) -- (1/2,1);

    \begin{scope}[thick]
        \draw[fill=gray!10] (0.0,0.0) rectangle (1/3,1/4);
        \draw[fill=gray!10] (0.0,1/2) rectangle (1/3,1/2+1/10);

        \draw[fill=gray!10] (1/2,1/4-1/7) rectangle (1,1/4+1/3-1/7);
        \draw[fill=gray!10] (1/2,3/4) rectangle (1,1);
    \end{scope}
\end{tikzpicture}

%% file: figures/carpet_attractor/fig.tex
\begin{tikzpicture}[scale=5]
    \begin{scope}[dotted]
        \draw[fill=gray!5] (0.0,0.0) rectangle (1/3,1/4);
        \draw[fill=gray!5] (0.0,1/2) rectangle (1/3,1/2+1/10);

        \draw[fill=gray!5] (1/2,1/4-1/7) rectangle (1,1/4+1/3-1/7);
        \draw[fill=gray!5] (1/2,3/4) rectangle (1,1);
    \end{scope}

    \begin{scope}[font=\tiny]
        \node[below] at (0,0) {$0$};
        \node[below] at (1,0) {$1$};
        \node[left] at (0,0) {$0$};
        \node[left] at (0,1) {$1$};
    \end{scope}
    \draw (0,0) rectangle (1,1);

    \draw[dashed] (1/3,0) -- (1/3,1);
    \draw[dashed] (1/2,0) -- (1/2,1);

    \node (gl) at (1/2,1/2) {\includegraphics[width=5cm]{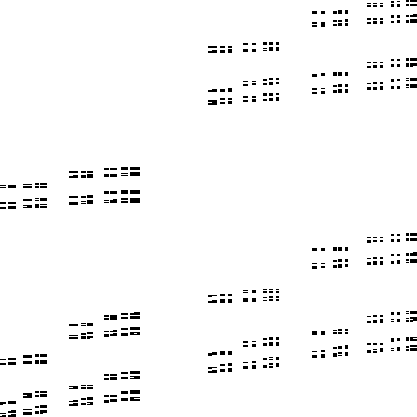}};
\end{tikzpicture}

%% file: figures/piecewise/fig.tex
\begin{tikzpicture}[xscale=40,yscale=10]
    \coordinate (t0) at (0.118585, -0.480531);
    \coordinate (t1) at (0.166828, -0.283151);
    \coordinate (t2) at (0.430605, 0);

    \coordinate (i0) at (0.118585, -0.53);
    \coordinate (i1) at (0.166828, -0.53);
    \coordinate (i2) at (0.430605, -0.53);

    \draw (i0) -- node[fill=white]{$I_1$} (i1) -- node[fill=white]{$I_2$} (i2);

    \coordinate (sl0) at (0.0873834, -0.684913);
    \coordinate (sr0) at (0.191945, 0);

    \coordinate (sl1) at (0.0873834, -0.475888);
    \coordinate (sr1) at (0.283541, 0);

    \coordinate (sl2) at (0.0873834, -0.414791);
    \coordinate (sr2) at (0.33771, 0);

    \coordinate (sl3) at (0.0873834, -0.343306);
    \coordinate (sr3) at (0.430605, 0);

    \draw[dashed,<->] (0.10,0) node[left]{$0$} -- (0.44,0);

    \draw[dashed] (i0) node[below=0.5cm]{$t_{\min}$} -- (t0);
    \draw[dashed] (i1) node[below=0.5cm]{$t_1$} -- (t1);
    \draw[dashed] (i2) node[below=0.5cm]{$t_{\max}$} -- (t2);
    \draw[thick, localorange] plot file {figures/piecewise/transform_1.txt};
    \draw[thick, localgreen, smooth] plot file {figures/piecewise/transform_2.txt};

    \draw[thick, draw=gray, dotted] (sl0) -- (sr0) node[anchor=south east,fill=white,rotate=58.59]{$\phi(\theta_{1,\min})$};
    \draw[thick,draw=gray,dotted] (sl1) -- (sr1) node[anchor=south east,fill=white,rotate=31.24]{$\phi(\theta_{1,\max})$};
    \draw[thick,draw=gray,dotted] (sl2) -- (sr2) node[anchor=south east,fill=white,rotate=22.5]{$\phi(\theta_{2,\min})$};
    \draw[thick,draw=gray,dotted] (sl3) -- (sr3) node[anchor=south east,fill=white,rotate=14.04]{$\phi(\theta_{2,\max})$};

    \node[vtx] (in0) at (i0){};
    \node[vtx] (in1) at (i1){};
    \node[vtx] (in2) at (i2){};

    \node[vtx] (tn0) at (t0){};
    \node[vtx] (tn1) at (t1){};
    \node[vtx] (tn2) at (t2){};
\end{tikzpicture}

%% file: figures/two_bump_spectrum/fig.tex
\begin{tikzpicture}
    \begin{axis}[
        ymax=1.057,
        xmax=1.05,
        xtick={0.001,0.999},
        xticklabels={$0$,$1$},
        xtick pos=bottom,
        ytick={1.0399520216923197,1.05528715},
        yticklabels={$\dimB K$, $\dimA K$},
        width=0.95\textwidth,
        height=6cm,
        xlabel={$\theta$},
        axis x line=center,
        axis y line=center
        ]
        \addplot[thick,black] table [x=x, y=y, col sep=comma] {figures/two_bump_spectrum/transition_1.txt};
        \addplot[thick,black] table [x=x, y=y, col sep=comma] {figures/two_bump_spectrum/transition_2.txt};
        \addplot[thick,black] table [x=x, y=y, col sep=comma] {figures/two_bump_spectrum/transition_3.txt};
        \addplot[thick,localorange] table [x=x, y=y, col sep=comma] {figures/two_bump_spectrum/conjugate_1.txt};
        \addplot[thick,localgreen] table [x=x, y=y, col sep=comma] {figures/two_bump_spectrum/conjugate_2.txt};

        \draw[localblue] (0.19,1.053) -- (0.19,1.0555) -- (0.28,1.0555) -- (0.28,1.053) -- cycle;
        \node[below] (ylabel) at (0.9,1.05428715) {$\dimAs\theta K$};
    \end{axis}
\end{tikzpicture}

%% file: figures/two_bump_spectrum_zoomed/fig.tex
\begin{tikzpicture}
    \begin{axis}[
        no markers,
        axis line style={draw=none},
        xtick={0.207523, 0.220413, 0.246873, 0.253354},
        xtick style={draw=none},
        ytick style={draw=none},
        xtick pos=bottom,
        xticklabels={$\theta_{1,\min}$,$\theta_{1,\max}$,$\theta_{2,\min}$,$\theta_{2,\max}$},
        yticklabels=\empty,
        width=16cm,
        height=6cm,
        axis x line=center %
        ]
        \draw[dashed] (0.21,1.05528715) node[left]{$\dimA K$} -- (0.253354,1.05528715);

        \addplot[thick,black,<-,>=stealth] table [x=x, y=y, col sep=comma] {figures/two_bump_spectrum_zoomed/transition_1.txt};
        \addplot[thick,black] table [x=x, y=y, col sep=comma] {figures/two_bump_spectrum_zoomed/transition_2.txt};
        \draw[thick,black,->,>=stealth] (0.253354,1.05528715) -- (0.268,1.05528715);
        \addplot[thick,localorange] table [x=x, y=y, col sep=comma] {figures/two_bump_spectrum_zoomed/conjugate_1.txt};
        \addplot[thick,localgreen] table [x=x, y=y, col sep=comma] {figures/two_bump_spectrum_zoomed/conjugate_2.txt};

        \draw[dashed] (0.207523, 1.0531) -- (0.207523,1.05395);
        \draw[dashed] (0.220413, 1.0531) -- (0.220413,1.05464);
        \draw[dashed] (0.246873, 1.0531) -- (0.246873,1.05521);
        \draw[dashed] (0.253354, 1.0531) -- (0.253354,1.05529);
    \end{axis}
    \draw[localblue] (current bounding box.north east) -- (current bounding box.north west) -- (current bounding box.south west) -- (current bounding box.south east) -- cycle;
\end{tikzpicture}

%% file: figures/homog_transform/fig.tex
\begin{tikzpicture}[xscale=500,yscale=60]
    \begin{scope}[thick, domain=0.274976:0.30103]
        \draw[localorange] plot (\x,2.71340271*\x-0.7990164);
        \draw[localgreen] plot (\x,1.92926933*\x-0.57571664);
        \draw[localblue] plot (\x,1.04318002*\x-0.31402848);
    \end{scope}

    \coordinate (b0) at (0.274976,0);
    \coordinate (b1) at (0.284773,0);
    \coordinate (b2) at (0.295329,0);
    \coordinate (b3) at (0.30103,0);

    \coordinate (t0) at (0.274976,-0.052895);
    \coordinate (t1) at (0.284773,-0.0263135);
    \coordinate (t2) at (0.295329,-0.00594682);
    \coordinate (t3) at (0.30103, 0);

    \draw[dashed] (b0) node[above,fill=white]{$t_0$} -- (t0);
    \draw[dashed] (b1) node[above,fill=white]{$t_1$} -- (t1);
    \draw[dashed] (b2) node[above,fill=white]{$t_2$} -- (t2);
    \draw[dashed] (b3) node[above,fill=white]{$t_3$} -- (t3);

    \node[vtx] (n0) at (b0){};
    \node[vtx] (n1) at (b1){};
    \node[vtx] (n2) at (b2){};
    \node[vtx] (n3) at (b3){};

    \node[vtx] (nt0) at (t0){};
    \node[vtx] (nt1) at (t1){};
    \node[vtx] (nt2) at (t2){};

    \draw (b0) -- (b3);
\end{tikzpicture}

%% file: figures/homog_spectrum/fig.tex
\def\rhomax{0.63145869}
\def\baseline{1.1758}
\begin{tikzpicture}[xscale=12,yscale=150]
    \begin{scope}[thick]
        \draw[domain=0:0.368541] plot (\x, 0.900824+0.274976+0.052895*\rhomax*\x/(1-\x);
        \draw[domain=0.368541:0.48,dotted] plot (\x, 0.900824+0.274976+0.052895*\rhomax*\x/(1-\x);

        \draw[domain=0:0.368541,dotted] plot (\x, 0.900824+0.284773+0.0263135*\rhomax*\x/(1-\x);
        \draw[domain=0.368541:0.450805] plot (\x, 0.900824+0.284773+0.0263135*\rhomax*\x/(1-\x);
        \draw[domain=0.450805:0.56,dotted] plot (\x, 0.900824+0.284773+0.0263135*\rhomax*\x/(1-\x);

        \draw[domain=0.602873:0.74,dotted] plot (\x, 0.900824+0.295329+0.00594682*\rhomax*\x/(1-\x);
        \draw[domain=0.450805:0.602873] plot (\x, 0.900824+0.295329+0.00594682*\rhomax*\x/(1-\x);
        \draw[domain=0:0.450805,dotted] plot (\x, 0.900824+0.295329+0.00594682*\rhomax*\x/(1-\x);

        \draw[dotted] (0,0.900824 + 0.30103) node[left]{$\dimA K$}-- (0.602873,0.900824 + 0.30103);
        \draw (0.602873,0.900824 + 0.30103) -- node[above]{$\dimAs\theta K$} (1,0.900824 + 0.30103);
    \end{scope}

    \draw[thick,->] (0,\baseline) node[left]{$\dimB K$} -- (1.05,\baseline);
    \draw (0,\baseline+0.0008) -- (0,\baseline-0.0008) node[below]{$0$};
    \draw (1,\baseline+0.0008) -- (1,\baseline-0.0008) node[below]{$1$};

    \draw[dashed] (0.368541,\baseline-0.0008) node[below]{$\theta_1$} -- (0.368541,1.19529);
    \draw[dashed] (0.450805,\baseline-0.0008) node[below]{$\theta_2$} -- (0.450805,1.19924);
    \draw[dashed] (0.602873,\baseline-0.0008) node[below]{$\theta_3$} -- (0.602873,1.20185);

    \node[vtx, localorange] (s1) at (0.368541,1.19529) {};
    \node[vtx, localgreen] (s2) at (0.450805,1.19924) {};
    \node[vtx, localblue] (s3) at (0.602873,1.20185) {};
\end{tikzpicture}

%% file: figures/two_int/fig.tex
\def\baseline{0.9910}
\begin{tikzpicture}
    \begin{axis}[
        no markers,
        xtick style={draw=none},
        ytick style={draw=none},
        xtick pos=bottom,
        xtick={0.268985, 0.271996, 0.277686, 0.280779},
        xticklabels={$\theta_{1,\max}$, $\theta_{2,\min}$, $\theta_{2,\max}$, $\theta_{3,\min}$},
        yticklabels=\empty,
        width=12cm,
        height=10cm,
        xlabel={$\theta$},
        ylabel={$\dimAs\theta K$},
        axis lines = middle,
        xmin=0.2673
        ]
        \addplot[thick,localorange] table [x=x, y=y, col sep=comma] {figures/two_int/conjugate_1.txt};
        \addplot[thick,localgreen] table [x=x, y=y, col sep=comma] {figures/two_int/conjugate_2.txt};
        \addplot[thick,black] table [x=x, y=y, col sep=comma] {figures/two_int/transition_1.txt};
        \addplot[thick,black] table [x=x, y=y, col sep=comma] {figures/two_int/transition_2.txt};

        \draw[dashed] (0.268985, \baseline) -- (0.268985,0.99141);
        \draw[dashed] (0.271996, \baseline) -- (0.271996,0.991633);
        \draw[dashed] (0.277686, \baseline) -- (0.277686,0.991911);
        \draw[dashed] (0.280779, \baseline) -- (0.280779,0.991983);
    \end{axis}
\end{tikzpicture}

%% file: figures/convex_transform/fig.tex
\begin{tikzpicture}[xscale=35,yscale=6]
    \draw[thick, localorange, smooth] plot file {figures/convex_transform/transform_1.txt};

    \coordinate (i0) at (0.0843522,0);
    \coordinate (t0) at (0.0843522,-0.63224);
    \coordinate (i1) at (0.417375,0);
    \coordinate (t1) at (i1);

    \draw[dashed] (i0) node[above=0.2cm,fill=white]{$t_0$} -- (t0);
    \path (i1) node[above=0.2cm, fill=white]{$t_1$} -- (t1);

    \draw (i0) -- (i1);

    \node[vtx] (n0) at (i0) {};
    \node[vtx] (nt0) at (t0) {};
    \node[vtx] (n1) at (i1) {};
\end{tikzpicture}

%% file: figures/convex_spectrum/fig.tex
\begin{tikzpicture}[xscale=12,yscale=12]
    \draw[thick, localorange, smooth, dotted] plot file {figures/convex_spectrum/conjugate_dotted_1.txt};
    \draw[thick, localorange, smooth] plot file {figures/convex_spectrum/conjugate_1.txt};
    \draw[thick, smooth] plot file {figures/convex_spectrum/transition_1.txt};
    \draw[thick] (0.500364,1.41737) -- (1,1.41737);

    \draw[dashed] (0.0839257, 1.08435) node[below]{$\theta_{1,\min}$} -- (0.0839257,1.1404);
    \draw[dashed] (0.500364, 1.08435) node[below]{$\theta_{1,\max}$} -- (0.500364,1.41737);

    \draw[->] (-0.02, 1.08435) -- (1.02, 1.08435) node[above left]{$\theta$};
    \draw[->] (0, 1.08435-0.02) node[below]{$0$} -- (0, 1.41737+0.02) node[below right]{$\dimAs\theta K$};
    \draw (1, 1.08435-0.01) node[below]{$1$} -- (1, 1.08435+0.01);
\end{tikzpicture}